\documentclass[10pt,a4paper]{article}
\usepackage{amsmath,amssymb,amsfonts,amsthm}
\usepackage{float,graphicx,color}
\usepackage[all]{xy}
\usepackage{tikz}

\newcommand{\rmd}{\mathrm{d}}

\newcommand{\rmD}{\mathrm{D}}
\newcommand{\rmH}{\mathrm{H}}
\newcommand{\rmL}{\mathrm{L}}

\newcommand{\rmW}{\mathrm{W}}

\newcommand{\bbA}{\mathbb{A}}
\newcommand{\bbB}{\mathbb{B}}
\newcommand{\bbC}{\mathbb{C}}
\newcommand{\bbE}{\mathbb{E}}

\newcommand{\bbN}{\mathbb{N}}
\newcommand{\bbP}{\mathbb{P}}

\newcommand{\bbR}{\mathbb{R}}
\newcommand{\bbS}{\mathbb{S}}

\newcommand{\frH}{\mathfrak{H}}

\newcommand{\calC}{\mathcal{C}}
\newcommand{\calD}{\mathcal{D}}
\newcommand{\calE}{\mathcal{E}}
\newcommand{\calF}{\mathcal{F}}

\newcommand{\calI}{\mathcal{I}}
\newcommand{\calL}{\mathcal{L}}
\newcommand{\calM}{\mathcal{M}}

\newcommand{\calS}{\mathcal{S}}
\newcommand{\calT}{\mathcal{T}}
\newcommand{\calU}{\mathcal{U}}
\newcommand{\calV}{\mathcal{V}}

\newcommand{\calY}{\mathcal{Y}}
\newcommand{\calZ}{\mathcal{Z}}

\newcommand{\llp}{\vert \! \lfloor}
\newcommand{\rrp}{\rfloor \! \vert}

\newcommand{\orient}{o} 
\newcommand{\hodge}{\frH}

\newcommand{\ds}{\displaystyle}
\newcommand{\ts}{\textstyle}
\newcommand{\bs}{{\scriptscriptstyle \bullet}}

\newcommand{\cleq}{\preceq}

\newcommand{\ceq}{\simeq}

\newcommand{\rmins}[1]{\quad \textrm{#1} \quad}

\newcommand{\poly}{\bbP}    
\newcommand{\alt}{\bbA}       

\DeclareMathOperator{\Div}{div}

\DeclareMathOperator{\grad}{grad}
\DeclareMathOperator{\curl}{curl}

\renewcommand{\div}{\Div}

\newcommand{\func}[4]{\left \{\begin{array}{rl}
#1 & \to #2,\\
#3 & \mapsto #4. 
\end{array}\right .
}

\DeclareMathOperator{\id}{id}

\newtheorem{theorem}{Theorem}[section]
\newtheorem{lemma}[theorem]{Lemma}
\newtheorem{proposition}[theorem]{Proposition}
\newtheorem{corollary}[theorem]{Corollary}

\theoremstyle{definition}\newtheorem{definition}{Definition}
\newtheorem{remark}{Remark}[section]
\newtheorem{example}{Example}[section]

\begin{document}
\title{Finite element systems of differential forms}
\author{Snorre H. Christiansen}

\date{Appeared as section 5 in \cite{ChrMunOwr11}.}

\maketitle

\begin{abstract}
We develop the theory of mixed finite elements in terms of special inverse systems of complexes of differential forms, defined over cellular complexes. Inclusion of cells corresponds to pullback of forms. The theory covers for instance composite piecewise polynomial finite elements of variable order over polyhedral grids. Under natural algebraic and metric conditions, interpolators and smoothers are constructed, which commute with the exterior derivative and whose product is uniformly stable in Lebesgue spaces. As a consequence we obtain not only eigenpair approximation for the Hodge-Laplacian in mixed form, but also variants of Sobolev injections and translation estimates adapted to variational discretizations.
\end{abstract}

\section{Introduction}
Many partial differential equations (PDEs) can be naturally thought of as expressing that a certain field, say a scalar field or a vector field, is a critical point of a certain functional. Often this functional will be of the form:
\begin{equation}
\calS(u) = \int_S \calL(x, u(x), \nabla u (x), \ldots) \rmd x. 
\end{equation}
Here $S$ is a domain in $\bbR^d$, $u$ is a section of a vector bundle over $S$ and $\calL$ is a Lagrangian involving $x\in S$ as well as values of $u$ and its derivatives. In this context the functional $\calS$ is called the action. Criticality of the action can be written:
\begin{equation}
\forall u' \quad \rmD \calS (u) u' = 0. 
\end{equation}
More generally, a PDE written in the form $\calF(u) = 0$ can be given a variational formulation:
\begin{equation}\label{eq:variational}
\forall u' \quad \langle \calF(u), u' \rangle = 0.
\end{equation}
The operator $\calF$ is continuous from a Banach space $X$ to another one $Y$, and $\langle \cdot, \cdot \rangle$  denotes a duality product on $Y \times Y'$. The unknown $u$ is sought in $X$ and the equation is tested with $u'$ in $Y'$. The duality product is most often a generalization of the $\rmL^2(S)$ scalar product, such as a duality between Sobolev spaces. In the case of least action principles it is of interest to find $X$ such that one can take $Y= X'$ (and $Y' = X$).

This point of view leads quite naturally to discretizations. One constructs a trial space $X_n \subseteq X$ and a test space $Y_n' \subseteq Y'$ and solve:
\begin{equation}\label{eq:variationaldisc}
u \in X_n,\  \forall u' \in Y'_n \quad \langle \calF(u), u' \rangle = 0.
\end{equation}
In theory one has sequences indexed by $n\in \bbN$ and establishes convergence properties of the method as $n \to \infty$. This discretization technique is called the Galerkin method when $X_n = Y'_n$ and the Petrov- Galerkin method when $X_n$ and $Y'_n$ are different.

The finite element (FE) method consists in constructing finite dimensional trial and test spaces adapted to classes of PDEs. Mixed FE spaces have been constructed  to behave well for the differential operators $\grad$, $\curl$ and $\div$. They provide a versatile tool-box for discretizing PDEs expressed in these terms. Thus Raviart-Thomas $\div$-conforming FEs \cite{RavTho77} are popular for PDEs in fluid dynamics and N\'ed\'elec's $\curl$-conforming FEs \cite{Ned80} have imposed themselves in electromagnetics. For reviews and references see \cite{BreFor91}\cite{RobTho91}\cite{Hip02}\cite{Mon03}.

These spaces fit the definition of a FE given by Ciarlet \cite{Cia78}. In particular they are equipped with unisolvent degrees of freedom, which determine inter element continuity and provide interpolation operators --  projections $I_n$ onto $X_n$ which are defined at least for smooth fields. In applications one often needs a pair of spaces $X_n^a \times X_n^b$, and these spaces should be compatible in the sense of satisfying a Brezzi inf-sup condition \cite{Bre74}. The spaces are linked by a differential operator $\rmd : X^a_n \to X^b_n$ (one among $\grad$, $\curl$, $\div$), and the proof of compatibility would follow from a commuting diagram:
\begin{equation}
\xymatrix{
X^a \ar[r]^\rmd \ar[d]^{I^a_n} & X^b\ar[d]^{I^b_n}  \\  
X^a_n \ar[r]^\rmd & X^b_n 
}
\end{equation}
A technical problem is that one would like the interpolators to be defined and bounded on the Banach spaces $X^a$, $X^b$, whereas the natural degrees of freedom are not bounded on these. For instance $X^a$ is usually of the form:
\begin{equation}
X^a = \{u \in \rmL^2(S)\otimes E^a \ : \rmd u \in \rmL^2(S)\otimes E^b \}, 
\end{equation}
for some finite dimensional fibers $E^a$, $E^b$, but this continuity is usually insufficient for degrees of freedom, such as line integrals, to be well-defined. The convergence of eigenvalue problems requires even more boundedness properties of the interpolation operators, something more akin to boundedness in $\rmL^2(S)$ (see \cite{ChrWin13IMA} for precise statements).

Mixed FE have been constructed to contain polynomials of arbitrarily high degree. Applying the Bramble Hilbert lemma to the interpolators gives high orders of approximation, and under stability conditions this gives high order of convergence for numerical methods. But here too the intention is hampered by the lack of continuity of interpolators.

As remarked by Bossavit \cite{Bos88}, lowest order mixed FE, when translated into the language of differential forms, correspond to constructs from differential topology called Whitney forms \cite{Wei52} \cite{Whi57}. The above mentioned differential operators can all be interpreted as instances of the exterior derivative. From this point of view it becomes natural to arrange the spaces in full sequences linked by operators forming commuting diagrams:
\begin{equation}
\xymatrix{
X^0  \ar[d]^{I^0_n} \ar[r]^\rmd & X^1\ar[d]^{I^1_n} \ar[r]^\rmd & \ldots \ar[d] \ar[r]^\rmd   & X^d\ar[d]^{I^d_n}  \\  
X^0_n \ar[r]^\rmd & X^1_n \ar[r]^\rmd & \ldots  \ar[r]^\rmd   & X^d_n. 
}
\end{equation}
An important property is that the interpolators should induce isomorphisms on cohomology groups. This is essentially De Rham's theorem, when the top row consists of spaces of smooth differential forms and the bottom row consists of Whitney forms. We retain that on domains homeomorphic to balls, the sequences of FE spaces should be exact. This alpplies in particular to single elements such as cubes and simplexes.

High order FE spaces of differential forms naturally generalizing Raviart-Thomas-N\'ed\'elec elements were presented by Hiptmair \cite{Hip99}. For a comprehensive review, encompassing also Brezzi-Douglas-Marini elements \cite{BreDouMar85}, see  \cite{ArnFalWin06}. Elements are constructed using the Koszul operator, or equivalently the Poincar\'e homotopy operator, to ensure local sequence exactness - that is, exactness of the discrete sequence attached to a single element. Mixed FE spaces have thus been defined for simplexes. Tensor product constructions yield spaces on Cartesian products thereof. 

Given spaces defined on some mesh it, would sometimes be useful to have spaces constructed on the dual mesh, matching the spaces on the primal mesh in some sense. Since the dual mesh of a simplicial mesh is not simplicial, this motivates the construction of FE spaces on meshes consisting of general polytopes. On these, it seems unlikely that good FE spaces can be constructed with only polynomials, one should at least allow for piecewise such. In some situations, like convection dominated fluid flow, stability requires some form of upwinding. This could be achieved in a Petrov-Galerkin method, by including upwinded basis functions in either the trial space $X_n$ or the test space $Y'_n$. For a recent review on this topic see \cite{Mor10}. This provides another motivation for constructing a framework for FE that includes non-polynomial functions.

Discretization of PDEs expressed in terms of $\grad$, $\curl$ and $\div$ on polyhedral meshes has long been pursued in the framework of mimetic finite differences, reviewed in \cite{BocHym06}. The convergence of such methods has been analysed in terms of related mixed finite elements \cite{BreLipSha05}. For connections with finite volume methods see \cite{DroEymGalHer10}. While recent developments tend to exhibit similarities and even equivalences between all these methods, the FE method, as it is understood here, could be singled out by its emphasis on defining fields inside cells, and ensuring continuity properties between them, in such a way that the main discrete differential operators acting on discrete fields are obtained symply by restriction of the continuous ones.

Ciarlet's definition of a FE does not capture the fact that mixed FE behave well with respect to restriction to faces of elements. From the opposite point of view, once this is realized, it becomes natural to prove properties of mixed FE by induction on the dimension of the cell. 

The goal is then to construct a framework for FE spaces of differential forms on cellular complexes accommodating arbitrary functions.  One requires stability of the ansatz spaces under restriction to subcells and under the exterior derivative. For the good properties of standard spaces to be preserved, one must impose additional conditions on the ansatz spaces, essentially surjectivity of the restriction from the cell to the boundary, and sequence exactness under the exterior derivative on each cell. As it turns out these conditions, which we refer to as compatibility, imply the existence of interpolation operators commuting with the exterior derivative. Degrees of freedom are not part of the definition of compatible finite element systems but are rather deduced, and it becomes more natural to compare various degrees of freedom for a given system. The local properties of surjectivity and sequence exactness also ensure global topological properties by the general methods of algebraic topology.

Interpolation operators still lack desired continuity properties.  But combining them with a smoothing technique yields commuting projection operators that are stable in $\rmL^2(S)$. Stable commuting projections were first proposed in \cite{Sch08}. Smoothing was achieved by taking averages over perturbations of the grid. Another construction using cut-off and smoothing by convolution on reference macroelements was introduced in \cite{Chr07NM}. While commutativity was lost, the lack of it was controlled by an auxilliary parameter. As remarked in \cite{ArnFalWin06}, for quasiuniform meshes one can simplify these constructions and use smoothing by convolution on the physical domain. If in \cite{Sch08} one can say that the nodes of the mesh are shaken independently, smoothing by convolution consists in shaking them in parallel. Much earlier, in \cite{DodPat76}, convergence for the eigenvalue problem for the Hodge-Laplacian, discretized with Whitney forms, was proved  using smoothing by the heat kernel. For scalar functions, smoothing by convolution in the FE method has been used at least as far back as \cite{Str72NM} \cite{Hil73}, but Cl\'ement interpolation \cite{Cle75} seems to have supplanted it. In \cite{ChrWin08} a space dependent smoothing operator, commuting with the exterior derivative, was introduced, allowing for general shape regular meshes. These constructions also require a commuting extension operator, extending differential forms outside the physical domain.

This paper is organized as follows. In section \ref{sec:interpol}, cellular complexes and the associated framework of finite element systems are introduced. Basic examples are included, as well as some constructions like tensor products. Section \ref{sec:interpol} serves to introduce degrees of freedom and interpolation operators on FE systems. In section \ref{sec:quasiint}, we  construct smoothers and extensions which commute with the exterior derivative and preserve polynomials locally. When combined with interpolators they yield $\rmL^q(S)$ stable commuting projections for scale invariant FE systems. In section \ref{sec:sob} we apply these constructions to prove discrete Sobolev injections and a translation estimate.

The framework of FE systems was implicit in \cite{Chr08M3AS} and made explicit in \cite{Chr09AWM} but we have improved some of the proofs. Upwinding in this context is new, as well as the discussion of interpolation and degrees of freedom. The construction of smoothers and extensions improves that of \cite{ChrWin08} by having the additional property of preserving polynomials locally, up to any given maximal degree. The analysis is also extended from $\rmL^2$ to $\rmL^q$ estimates, for all finite $q$. This is used in the proof of the Sobolev injection and translation estimate which are also new (improving \cite{ChrSch11} and \cite{KarKar11}).


\section{Finite element systems \label{sec:fesys}}

\paragraph{Cellular complexes.}

For any natural number $k\geq 1$, let $\bbB^k$ be the closed unit ball of $\bbR^k$ and $\bbS^{k-1}$ its boundary. For instance $\bbS^0= \{-1, 1\}$. We also put $\bbB^0 = \{0\}$.

Let $S$ denote a compact metric space. A $k$-dimensional \emph{cell} in $S$ is a closed subset $T$ of $M$ for which there is a Lipschitz bijection $\bbB^k \to T$ with a Lipschitz inverse. If a cell $T$ is both $k$- and $l$-dimensional then $k= l$. For $k \geq 1$, we denote by $\partial T$ its boundary, the image of $\bbS^{k-1}$ by the chosen bi-Lipschitz map. Different such maps give the same boundary. The interior of $T$ is by definition $T \setminus \partial T$ (it is open in $T$ but not necessarily in $S$). 
\begin{definition}
A \emph{cellular complex} is a pair $(S,\calT)$ where $S$ is a compact metric space and $\calT$ is a finite set of cells in $S$, such that the following conditions hold:
\begin{itemize}
\item Distinct cells in $\calT$ have disjoint interiors.
\item The boundary of any cell in $\calT$ is a union of cells in $\calT$.
\item The union of all cells in $\calT$ is $S$.
\end{itemize}
\end{definition}
In this situation we also say that $\calT$ is a cellular complex on $S$. One first remarks:
\begin{proposition} The intersection of two cells in $\calT$ is a union of cells in $\calT$.
\end{proposition}
\begin{proof}
Let $T, U$ be two cells in $\calT$ and suppose $x \in T \cap U$. Choose a cell $T'$ included in $T$ of minimal dimension such that $x \in T'$. Choose also a cell $U'$ included in $U$ of minimal dimension such that $x \in U'$. Suppose none of the cells $T'$ and $U'$ are points. Then $x$ belongs to the interiors of both $T'$ and $U'$, so that $T' = U'$. Therefore there exists a cell included in both $T$ and $U$ to which $x$ belongs. This conclusion also trivially holds if $T'$ or $U'$ is a point.
\end{proof}

In fact, if $(S,\calT)$ is a cellular complex, $S$ can be recovered from $\calT$ as follows:
\begin{proposition} When $(S,\calT)$ is a cellular complex:
\begin{itemize}
\item $S$ is recovered as a set from $\calT$ as its union.
\item The topology of $S$ is determined by the property that a subset $U$ of $S$ is closed iff for any cell $T$ in $\calT$, of dimension $k$, the image of $U \cap T$ under the chosen bi-Lipschitz map, is closed in $\bbB^k$.
\item A compatible metric on $S$ can be recovered from $\calT$, from metrics $d_T$ on each cell $T$ inherited from $\bbB^k$, for instance by ($ \max \emptyset = + \infty$):
\begin{equation}
d(x,y) = \min \{1, \ \max \{ d_T(x,y) \ : \ x,y \in T \rmins{and}  T \in \calT \}.
\end{equation}
\end{itemize}

\end{proposition}

A \emph{simplicial complex} is a cellular complex in which the intersection of any two cells is a cell (not just a union of cells) and such that the boundary of any cell is split into subcells in the same way as the boundary of a reference simplex is split into subsimplexes. The reference simplex $\Delta_k$ of dimension $k$ is:
\begin{equation}
\{ (x_0, \ldots , x_k) \in \bbR^{k + 1} \ : \  \sum_i x_i = 1  \rmins{and} \forall i \quad x_i \geq 0 \}.
\end{equation}
Its subsimplexes are parametrized by subsets $J \subseteq \{0, \cdots, k\}$ as the solution spaces:
\begin{equation}
\{x \in \Delta_k \ : \  \forall i \not \in J \quad x_i = 0 \}.
\end{equation}

For each $k\in \bbN$ the subset of $\calT$ of $k$-dimensional sets is denoted:
\begin{equation}
\calT^k= \{T \in \calT \ : \ \dim T = k\}.
\end{equation}  
The $k$-\emph{skeleton} of $\calT$ is the cellular complex consisting of cells of dimension at most $k$:
\begin{equation}
 \calT^{(k)} = \calT^0 \cup \cdots \cup \calT^k.
\end{equation}
The boundary $\partial T$ of any cell $T$ of $\calT$ can be naturally equipped with a cellular complex, namely: 
\begin{equation}
\{T'\in \calT \ : \ T' \subseteq T \rmins{and} T' \neq T \}. 
\end{equation}
We use the same notation for the boundary of a cell and the cellular complex it carries.

A \emph{refinement} of a cellular complex $\calT$ on $S$ is a cellular complex $\calT'$ on $S$ such that each element of $\calT$ is the union of elements of $\calT'$. We will be particularly interested in simplicial refinements of cellular complexes.

A cellular \emph{sub}complex of a cellular complex $\calT$ on $S$, is a cellular complex $\calT'$ on some closed part $S'$ of $S$ such that  $\calT' \subseteq \calT$. For instance if $T \in \calT$ is a cell, its subcells form a subcomplex of $\calT$, which we denote by $\tilde{T}$. We have seen that the boundary of any cell $T \in \calT$ can be equipped with a cellular complex which is a subcomplex of $\tilde{T}$.

Fix  a cellular complex $(S,\calT)$. In the following we suppose that for each $T \in \calT$ of dimension $\geq 1$, the manifold $T$ has been oriented. The \emph{relative orientation} of two cells $T$ and $T'$ in $\calT$, also called the \emph{incidence number}, is denoted $\orient(T,T')$ and defined as follows. For any edge $e\in \calT^1$ its vertexes are ordered, from say $\dot{e}$ to $\ddot{e}$. Define $\orient(e, \dot{e})= -1$ and $\orient(e, \ddot{e})=1$. Concerning higher dimensional cells, fix $k \geq 1$.
Given $T\in \calT^{k+1}$ and $T'\in \calT^{k}$ such that $T' \subseteq T$ we define  $\orient(T,T')=1$ if $T'$ is outward oriented compared with $T$ and $\orient(T,T')=-1$ if it is inward oriented. For all $T, T'\in \calT$ not covered by these definitions we put $\orient (T,T') = 0$.

For each $k$, let $\calC^k(\calT)$ denote the set of maps $c:\calT^k \to \bbR$. Such maps associate a real number with each $k$-dimensional cell and are called $k$-\emph{cochains}. The \emph{coboundary} operator $\delta : \calC^k(\calT) \to \calC^{k+1}(\calT)$ is defined by:
\begin{equation}
(\delta c)_T= \sum_{T'\in \calT^k} \orient(T,T')c_{T'}.
\end{equation}
The space of $k$-cochains has a canonical basis indexed by $\calT^k$. The coboundary operator is the operator whose canonical matrix is the incidence matrix $\orient$, indexed by $\calT^{k+1} \times \calT^k$. We remark that the coefficients in the sum can be non-zero only when $T'\in \partial T \cap \calT^k$.

\begin{lemma}
We have that $\delta \delta = 0$ as a map $\calC^k(\calT) \to \calC^{k+2}(\calT)$.
\end{lemma}
\begin{proof}See e.g. Lemma 3.6 in \cite{Chr09AWM}.
\end{proof}
In other words the family $\calC^\bs(\calT)$ is a complex, called the cochain complex and represented by:
\begin{equation*}
0 \to \calC^0(\calT) \to \calC^1(\calT) \to \calC^2(\calT) \to \cdots
\end{equation*}

When $S$ is a smooth manifold we denote by $\Omega^k(S)$ the space of smooth differential $k$-forms on $S$. Differential forms can be mapped to cochains as follows.
Let $S$ be a manifold and $\calT$ a cellular complex on $S$. 
For each $k$ we denote by $\rho^k : \Omega^k(S) \to \calC^k(\calT)$ the De Rham map,  which is defined by:
\begin{equation}
\rho^k: u \mapsto (\int_T u )_{T\in \calT^k}.
\end{equation}

\begin{proposition}
For each $k$ the following diagram commutes:
\begin{equation}
\xymatrix{
\Omega^k(S) \ar[r]^\rmd \ar[d]^{\rho^k} & \Omega^{k+1}(S) \ar[d]^{\rho^{k+1}}\\
\calC^k(\calT) \ar[r]^{\delta} & \calC^{k+1}(\calT)
}
\end{equation}
\end{proposition}
\begin{proof}
This is an application of Stokes theorem.
\end{proof}

Suppose $\calT$ is a cellular complex equipped with an orientation (of the cells) and $\calT'$ is a cellular refinement also equipped with an orientation (for instance the same complex but with different orientations). For each cell $T\in \calT^k$ and each $T'\in \calT^{\prime k}$ define $\iota(T,T') = \pm 1$ if $T' \subseteq T$ and they have the same/different orientation, and $\iota(T,T') = 0$ in all other cases. 
\begin{proposition}\label{prop:refcom}
The map $\iota : \calC^\bs(\calT') \to \calC^\bs(\calT)$ defined by:
\begin{equation}
(\iota u)_T= \sum_{T'\in \calT'}\iota(T,T') u_{T'},
\end{equation}
is a morphism of complexes, meaning that $\iota$ and $\delta$ commute.
\end{proposition}
\begin{proof}
See e.g. Proposition 3.5 in \cite{Chr09AWM}.
\end{proof}
We also remark that for a manifold $S$ we have:
\begin{proposition}
The following diagram commutes:
\begin{equation}
\xymatrix{
               & \Omega(S) \ar[dl]_\rho \ar[dr]^\rho\\
\calC(\calT')\ar[rr]^\iota & & \calC(\calT)
}
\end{equation}
\end{proposition}


\paragraph{Element systems}
For any cell $T$, we denote by $\Omega^k_{s,q}(T)$ the space of differential $k$-forms on $T$ with $\rmW^{s,q}(T)$ Sobolev regularity, and put $\Omega^k_q(T)= \Omega^k_{0,q}(T)$.
Fix $q \in [1, +\infty[$ and define:
\begin{equation}\label{eq:xdef}
X^k(T)= \{u \in \Omega^k_q(T)\  :  \ \rmd u \in \Omega^k_q(T) \}. 
\end{equation}
When $i:T' \to T$ is an inclusion of cells and $u$ is a smooth enough form on $T$ we denote by $u|_{T'} = i^\star u$ the pullback of $u$ to $T'$. Thus we restrict to the subcell and forget about the action of $u$ on vectors not tangent to it. In the topology (\ref{eq:xdef}), restrictions to subcells $T'$ of codimension one are well-defined, for instance as elements of $\Omega^k_{-1,q}(T')$. When $T$ is a cell in a given cellular complex $\calT$ we may therefore set:
\begin{equation}
 \hat X^k(T) = \{u \in X^k(T) \ : \forall T' \in \calT  \quad T' \subseteq T \Rightarrow u|_{T'} \in X^k(T') \}.
\end{equation}

\begin{definition}
Suppose $\calT$ is a cellular complex. For each $k\in \bbN$ and each $T \in \calT$ we suppose we are given a space $A^k(T) \subseteq \hat X^k(T) $ called a differential $k$-element on $T$.  We suppose that the exterior derivative induces maps $\rmd: A^k(T) \to A^{k+1}(T)$ and that if $i: T' \subseteq T$ is an inclusion of cells, pullback induces a map $i^\star: A^k(T) \to A^k(T')$. Such a family of elements is called an \emph{element system}. 
\end{definition}

A differential element is said to be finite if it is finite dimensional. A finite element system is an element system in which all the elements are finite. 

\begin{example}
The spaces $\hat X^\bs(\bs)$ themselves define an element system.  It is far from finite. 
\end{example}

\begin{example}
Let $U$ be an open subset of a vector space $V$. We denote by $\poly_p(U)$ the space of real polynomials of degree at most $p$ on $U$. For $k \geq 1$ the space of alternating maps $V^k \to \bbR$ is denoted $\alt^k(V)$. The space of differential $k$-forms on $U$, which are polynomial of degree at most $p$, is denoted $\poly\alt^k_p(U)$. We identify:
\begin{equation}
\poly\alt^k_p(U)=\poly_p(U) \otimes \alt^k(V) \rmins{and} \poly\alt^0_p(U) = \poly_p(U). 
\end{equation}
Choose a cellular complex where all cells are flat. Choose a function $\pi: \calT \times \bbN \to \bbN$ and define: 
\begin{equation}
A^k(T) = \poly\alt^k_{\pi(T,k)}(T). 
\end{equation}
One gets a finite element system when the following conditions are satisfied:
\begin{equation}
T' \subseteq T \Rightarrow \pi(T',k) \geq \pi(T,k) \rmins{and} \pi(T,k+1) \geq \pi(T,k) -1.
\end{equation}
\end{example}

\begin{example}
Denote the Koszul operator on vector spaces by $\kappa$. It is the contraction of differential forms by the identity, considered as a vector field:
\begin{equation}
(\kappa u)_x(\xi_1, \ldots , \xi_k) = u_x(x,\xi_1, \ldots , \xi_k).
\end{equation}
Alternatively one can use the Poincar\'e operator associated with the canonical homotopy from the identity to the null-map.  Let $\calT$ be a simplicial complex. Define, for non-zero $p\in \bbN$:
\begin{align}
\Lambda^k_p(T) &= \{u \in \poly\alt^k_p(T) \ : \  \kappa u \in \poly\alt^{k-1}_p(T) \} = \poly\alt^k_{p-1}(T) + \kappa \poly\alt^{k+1}_{p-1}(T).
\end{align}
For fixed $p$ we call this the \emph{trimmed} polynomial finite element system of order $p$. The case $p=1$ corresponds to constructs in \cite{Wei52}\cite{Whi57}. Arbitrary order elements were introduced in \cite{Ned80} for vector fields in $\bbR^3$. In \cite{Hip99} these spaces were extended to differential forms. The correspondence between lowest order mixed finite elements and Whitney forms was pointed out in \cite{Bos88}. See \cite{ArnFalWin06} for a comprehensive review. It was usual to start the indexing at $p=0$ but, as remarked in the preprint of \cite{Chr07NM}, the advantage of letting the lowest order be $p=1$, is that the wedge product induces maps:
\begin{equation}
\wedge : \Lambda^{k_0}_{p_0}(T) \times \Lambda^{k_1}_{p_1}(T) \to \Lambda^{k_0 + k_1}_{p_0 + p_1}(T).
\end{equation}
See also \cite{ArnFalWin06} p. 34. In words, the wedge product respects the grading in $k$ and the filtering in $p$. This observation was useful in the implementation of a scheme for the Yang-Mills equation \cite{ChrWin06}.
\end{example}

The first example of a finite element system yields quite useless Galerkin spaces in general, whereas the second one yields good ones. We shall elaborate on this in what follows, starting by defining what the Galerkin space associated with a finite element system is.

For any subcomplex $\calT'$ of $\calT$ we define $A^k(\calT')$ as follows :
\begin{equation}
A^k(\calT')= \{u \in \bigoplus_{T \in \calT'} A^k(T) \ : \  \forall T, T' \in \calT \quad  T' \subseteq T \Rightarrow u_{T}|_{T'} = u_{T'}\}. 
\end{equation}
Elements of $A^k(\calT')$ may be regarded as differential forms defined piecewise, which are continuous across interfaces between cells, in the sense of equal pullbacks. For a cell $T$ its collection of subcells is the cellular complex $\tilde{T}$. Applied to this case, the above definition gives a space canonically isomorphic to $A^k(T)$.  We can identify $A^k(\tilde T) = A^k(T)$.

A FE system over a cellular complex is an \emph{inverse system} of complexes: to an inclusion of cells corresponds the restriction operator. The space $A^\bs(\calT')$ defined above is an \emph{inverse limit} of this system and is determined by this property up to unique isomorphism.

Of particular importance is the application of the above construction to the boundary $\partial T$ of a cell $T$, considered as a cellular complex consisting of all subcells of $T$ except $T$ itself. Considering $\partial T$ as a cellular complex (not only a subset of $T$) we denote the constructed space by $A^k(\partial T)$.
If $i: \partial T \to T$ denotes the inclusion map, the pullback by $i$ defines a map $i^\star : A^k(T) \to A^k(\partial T)$ which we denote by $\partial$ and call restriction.

\emph{Convention:} in the following, the arrows starting or ending in $0$ are the only possible ones. Arrows starting in $\bbR$ are, unless otherwise specified, the maps taking a value to the corresponding constant function. Arrows ending in $\bbR$ are integration of forms of maximal degree. Other unspecified arrows are instances of the exterior derivative.

Consider now the following conditions on an element system $A$ on a cellular complex $\calT$: 
\begin{itemize}
\item \emph{Extensions.} For each $T\in \calT$ and $k\in \bbN$, restriction $\partial : A^k(T) \to A^k(\partial T)$ is onto.
\item \emph{Exactness.} The following sequence is exact for each $T$:
\begin{equation}\label{eq:coh}
0 \to \bbR \to A^0(T) \to A^1(T) \to \cdots \to A^{\dim T}(T) \to 0.
\end{equation}
\end{itemize}
The first condition can be written symbolically $\partial A^k(T) = A^k(\partial T)$.

\begin{definition}
We will say that an element system \emph{admits extensions} if the first property holds, is \emph{locally exact} if the second condition holds and is \emph{compatible} if both hold.
\end{definition}

Given a finite element system $A$, we say that its points carry reals, if for each \emph{point} $T\in \calT^0$, $A^0(T)$ contains the constant maps $T \to \bbR$ (so that $A^0(T)= \bbR^T \approx \bbR$).

\begin{proposition}\label{prop:one}
If $A$ admits extensions and its points carry reals, then for each cell $T$, $A^{\dim T}(T)$ contains a form with non-zero integral. 
\end{proposition}
\begin{proof}
By induction on the dimension of the cell, using Stokes' theorem.
\end{proof}

\emph{Notation:} We denote by  $A_0^k(T)$ be the kernel of $\partial: A^k(T) \to A^k(\partial T)$.

\begin{proposition}\label{prop:dima} We have:
\begin{equation}
\dim A^k(\calT) \leq  \sum_{T\in \calT} \dim A^k_0(T),
\end{equation}
with equality when the finite element system admits extensions.
\end{proposition}
\begin{proof} For a given $m \geq 0$, let $\calT^{(m)}$ be the $m$-skeleton of $\calT$. We have a sequence:
\begin{equation}
0 \to \bigoplus_{T \in \calT^{m}} A^k_0(T) \to A^k(\calT^{(m)})\to A^k(\calT^{(m-1)}) \to 0 . 
\end{equation}
The second arrow is bijective onto the kernel of the third. If all cells of dimension $m$ admit extensions the whole sequence is exact. The proposition follows from applying these remarks for all $m$.
\end{proof}

Combining Propositions \ref{prop:one} and \ref{prop:dima} we get:
\begin{corollary}
If $A$ admits extensions and its points carry reals, then $A^{k}(T)$ has dimension at least the number of $k$-dimensional subcells of  $T$.
\end{corollary}

\begin{proposition}\label{prop:gluecoh}
When the element system is compatible, the De Rham map $\rho^\bs : A^\bs(\calT) \to \calC^\bs(\calT)$ induces isomorphisms in cohomology.
\end{proposition}
\begin{proof}
By induction on the dimension of $\calT$. In dimension $0$ it is clear. Suppose now $m\geq 1$ and that we have proved the theorem when $\dim \calT < m$. Suppose that $\calT$ has dimension $m$. 

Remark that the De Rham map gives isomorphisms in cohomology on cells:
\begin{equation}
A^\bs(T) \to \calC^\bs(\tilde{T}),
\end{equation}
since both complexes are acyclic (?).

Denote by $\calU$ the $(m-1)$-skeleton of  $\calT$. Consider the diagram:
\begin{equation}
 \xymatrix{
0 \ar[r] & A^\bs(\calT^{}) \ar[r]  \ar[d]& A^\bs (\calU) \ds \bigoplus_{T \in \calT^{m}} A^\bs(T)   \ar[r] \ar[d]& \ds \bigoplus_{T \in \calT^{m}} A^\bs(\partial T)  \ar[r] \ar[d]& 0\\
0 \ar[r] & \calC^\bs(\calT^{}) \ar[r]  & \calC^\bs (\calU) \ds \bigoplus_{T \in \calT^{m}} \calC^\bs(\tilde{T})   \ar[r] & \ds \bigoplus_{T \in \calT^{m}} \calC^\bs(\partial T)  \ar[r] & 0
}
\end{equation}
The vertical maps are De Rham maps. The second horizontal arrow consists in restricting to the summands whereas the third one consists in restricting and comparing, as in the Mayer-Vietoris sequence.

Both rows are exact sequences of complexes, the diagram commutes, and the last two vertical arrows induce isomorphisms in cohomology by the induction hypothesis. Write the long exact sequences of cohomology groups associated with both rows (e.g. \cite{Chr09AWM} Theorem 3.1) and connect them with the induced morphisms.  Applying the five lemma (e.g. \cite{Chr09AWM} Lemma 3.2) gives the result for $\calT$.
\end{proof}

\begin{proposition}\label{prop:equiv}
For an element system with extensions the exactness of (\ref{eq:coh}) on each $T\in \calT$ is equivalent to the combination of the following two conditions:
\begin{itemize}
\item For each $T \in \calT$,  $A^0(T)$ contains the constant functions.
\item For each $T\in \calT$, the following sequence (with boundary condition) is exact:
\begin{equation} \label{eq:coh0}
0 \to A^0_0(T) \to A^1_0(T) \to \cdots \to A^{\dim T}(T) \to \bbR \to 0.
\end{equation}
\end{itemize}

 \end{proposition}
\begin{proof} When the extension property is satisfied on all cells, both versions (with and without boundary condition) of the cohomological condition guarantee that for each $T$, $A^{\dim T}(T)$ contains a form with non-zero integral, by Proposition \ref{prop:one}. 

When the spaces $A^0(T)$ all contain the constant functions, we may consider the following diagrams, where $T$ is a cell of dimension $m$ and $A^m_-(T)$ denotes the space of forms with zero integral:
\begin{equation}\label{eq:threerows}
\xymatrix{
0 \ar[r] & \bbR \ar[r]            & A^0 (\partial T)    \ar[r]  & \cdots \ar[r]          & A^{m-1}(\partial T) \ar[r]  & \bbR               \ar[r] &  0\\
0 \ar[r] & \bbR \ar[r]\ar[u]   & A^0(T)          \ar[r]\ar[u] & \cdots \ar[r]\ar[u] & A^{m-1}(T)     \ar[r] \ar[u]  & A^m(T) \ar[r]\ar[u] &  0\\
             &      0    \ar[r]\ar[u] & A^0_0(T)     \ar[r]\ar[u] & \cdots \ar[r]\ar[u] & A^{m-1}_0(T) \ar[r] \ar[u]  & A^m_{-}(T) \ar[r]\ar[u] &  0
}
\end{equation}
The columns, extended by $0$, are exact iff the extension property holds on $T$ and $A^m(T)$ contains a form of nonzero integral.
In this case, if one row is exact, the two other rows are either both exact or both inexact.

We now prove the stated equivalence.

(i) If compatibility holds then, in (\ref{eq:threerows}) the extended columns are exact, as well as the first and second row, so also the third. Hence (\ref{eq:coh0}) is exact.

(ii) Suppose now exactness of (\ref{eq:coh0}) holds for each $T$.  Choose $m\geq 1$ and suppose that we have proved exactness of (\ref{eq:coh}) for cells of dimension up to $m-1$. Let $T$ be a cell of dimension $m$.  In (\ref{eq:threerows}) the extended columns are exact. Apply Proposition \ref{prop:gluecoh} to the boundary of $T$ to get exactness of the first row. The third row is exact by hypothesis, and we deduce exactness of the second. The induction, whose initialization is trivial, completes.
\end{proof}

\begin{corollary} $\hat X$ is a compatible element system.
\end{corollary}



\paragraph{Tensor products.}

Suppose we have two manifolds $M$ and $N$, equipped with cellular complexes $\calU$ and $\calV$. We suppose we have differential elements $A^k(U)$ for $U \in \calU$ and $B^k(V)$ for $V \in \calV$, both forming systems as defined above.

Let $\calU \times \calV$ denote the product cellular complex on $M \times N$, whose cells are all those of the form $U \times V$ for $U \in \calU$ and $V \in \calV$. Recall that the tensor product of differential forms $u$ on $U$ and $v$ on $V$, is the form on $U \times V$ defined as the wedge product of their pullbacks by the respective canonical projections $p_U:U\times V \to U$ and $p_V:U \times V \to V$. In symbols we can write:
\begin{equation}
 u \otimes v = (p_U^\star u) \wedge (p_V^\star v).
\end{equation}
We equip $\calU \times \calV$ with elements:
\begin{equation}
C^\bs(U \times V) = A^\bs(U) \otimes B^\bs (V).
\end{equation}
Explicitly we put:
\begin{equation}
C^k(U \times V) = \bigoplus_l A^{l}(U) \otimes B^{k-l} (V).
\end{equation}
This defines an element system $C$, called the tensor product of $A$ and $B$.

\begin{proposition} \label{prop:czero} We have:
\begin{equation}
C^\bs_0(U \times V) = A^\bs_0(U) \otimes B^\bs_0(V).
\end{equation}
\end{proposition}
\begin{proof}
See \cite{Chr09AWM} Lemma 3.10.
\end{proof}

\begin{proposition}\label{prop:tensor}
When $A$ and $B$ admit extensions we have:
\begin{equation}
C^\bs(\calU \times \calV) = A^\bs(\calU) \otimes B^\bs(\calV).
\end{equation}
\end{proposition}
\begin{proof}
The right hand size is included in the left hand side. Moreover, by Proposition \ref{prop:dima}:
\begin{equation}
\dim C(\calU \times \calV)  \leq \ts \sum_{U,V} \dim C_0(U\times V),
\end{equation}
On the other hand Proposition \ref{prop:czero} gives:
\begin{align}
\ts \sum_{U,V} \dim C_0(U\times V) & =  \ts \sum_{U,V} \dim  A_0(U) \otimes B_0(V), \\
&=  \ts \sum_{U,V} \dim  A_0(U) \dim B_0(V),\\
&=  \ts \sum_U \dim  A_0(U) \ts \sum_V \dim B_0(V), \\
&= \dim A(\calU) \dim B(\calV).
\end{align}
This completes the proof.
\end{proof}

\begin{proposition}\label{prop:tensorext}
If $A$ and $B$ admit extensions, so does their tensor product.
\end{proposition}
\begin{proof}
Consider cells $U \in \calU$ and $V\in \calV$. Remark that:
\begin{align}
(\partial U \times V) \cup (U \times \partial V) & = \partial (U \times V),\\
(\partial U \times V) \cap (U \times \partial V) & = \partial U \times \partial V.
\end{align}
The Mayer-Vietoris principle gives an exact sequence:
\begin{equation}
0 \to C(\partial (U\times V)) \to C(\partial U \times V) \oplus C( U \times \partial V) \to C (\partial U \times \partial V) \to 0,
\end{equation}
where the second and third mappings are:
\begin{align}
w & \mapsto w|_{\partial U \times V} \ \oplus \ w|_{U \times \partial V},\\
u \oplus v & \mapsto u|_{\partial U \times \partial V} - v|_{\partial U \times \partial V}.
\end{align}
It follows that:
\begin{equation}
\dim C(\partial (U\times V))  = \dim C(\partial U \times V) + \dim C( U \times \partial V) -\dim C (\partial U \times \partial V).
\end{equation}
Applying Proposition \ref{prop:tensor} three times we get:
\begin{align}
\dim C(\partial (U\times V)) & = \dim A(\partial U) \dim B(V) + \dim A(U) \dim B(\partial V) - \\
& \quad  \quad \dim A(\partial U) \dim B(\partial V),\\
& = \dim A(U) \dim B(V) - \dim A_0(U) \dim B_0(V),\\
& = \dim C(U \times V) - \dim C_0(U \times V),\\
& = \dim \partial C(U \times V).
\end{align} 
Therefore:
\begin{equation}
\partial C(U \times V ) = C(\partial (U \times V)),
\end{equation}
as announced.
\end{proof}

\begin{proposition}
If $A$ and $B$ are locally exact, then so is their tensor product.
\end{proposition}
\begin{proof}
This follows from the Kunneth theorem (e.g. \cite{Chr09AWM} Theorem 3.2).
\end{proof}

\paragraph{Nesting.}

Suppose $\calT$ is a cellular complex and that $(\Xi, \preceq)$ is an ordered set. Suppose that, for each parameter $\xi \in \Xi$, a FE system $A[\xi]$ on $\calT$ has been chosen. We suppose that if $\xi \preceq \xi'$ then $A[\xi]^{k}(T) \subseteq A[\xi']^{k}(T)$. Choose now a parameter  function $\pi : \calT \to \Xi$ which is order preserving in the sense that if $T' \subseteq T$ then $\pi(T') \preceq \pi(T)$.
Define a FE system $A[\pi]$ by:
\begin{equation}
A[\pi]^k(T) = \{ u \in A[\pi(T)]^k(T) \ : \ \forall T' \in \calT \quad T' \subseteq  T \Rightarrow u |_{T'} \in A[\pi(T')]^{k} (T') \}.
\end{equation}

\begin{proposition}
If for each $\xi \in \Xi$, the system $A[\xi]$  is compatible, then the constructed system $A[\pi]$ too.
\end{proposition}
\begin{proof}
Remark that if $u \in A[\pi]^k(\partial T)$ then $u \in A[\xi(T)]^{k}(\partial T)$ so it can be extended to an element $u \in A[\xi(T)]^{k}(T)$. This element is in $A[\pi]^k(T)$. Thus $A[\pi]^k(T)$ admits extensions.

We also have $ A[\pi]^k_0(T) = A[\pi(T)]^{k}_0(T)$, which gives local exactness thanks to Proposition \ref{prop:equiv}.
\end{proof}

\begin{example} For a simplicial complex $\calT$, denote by $A[p]$ the trimmed finite element system of order $p$. Model variable order of approximation by a function $\pi : \calT \to \bbN^\ast$ such that when $T' \subseteq T$ we have $\pi(T') \leq \pi(T)$. The above construction defines a FE system $A[\pi]$, of the type used for $hp$-methods \cite{DemAl08}.  We will check later that the trimmed systems for fixed order $p$ are compatible, and then the above result gives compatibility of the variable order system defined by $\pi$.
\end{example}

\paragraph{Locally harmonic forms.}
Let $T$ be a cell where, for each $k$, $A^k(T)$ is equipped with a scalar product $a$. Orthogonality with respect to $a$ will be denoted $\perp$. We say that a $k$-form $u$ on $T$ is $A$--\emph{harmonic} if:
\begin{equation}\label{eq:aharmonic}
\rmd u \perp \rmd  A^k_0(T) \rmins{and} u \perp \rmd  A^{k-1}_0(T).
\end{equation}

One can for instance take $a$ to be the $\rmL^2$ scalar product on differential forms, associated with some Riemannian metric. Denote by $\rmd^\star$ the formal adjoint of $\rmd$ with respect to this scalar product. The continuous analogue of the above condition (\ref{eq:aharmonic}) is:
\begin{equation}\label{eq:l2harmonic}
\rmd^\star \rmd u = 0 \rmins{and} \rmd^\star u=0.
\end{equation}
From the other point of view, (\ref{eq:aharmonic}) is the Galerkin variant of (\ref{eq:l2harmonic}).

\begin{proposition}\label{prop:intext} Let $A$  be a finite element system where each  $A^k(T)$ is equipped with a scalar product $a$. Suppose $T$ is a cell such that (\ref{eq:coh0}) is exact. Put $m =\dim T$.
\begin{itemize}
\item For each  $\alpha \in \bbR$, there is a unique $A$-harmonic $u \in A^{m}(T)$ such that $\int_T u = \alpha $.
\item For $k<m$, any $u \in A^k(\partial T)$ admitting an extension in $A^k(T)$, has a unique $A$--harmonic extension in $A^k(T)$.
\end{itemize}
\end{proposition}
\begin{proof}
This follows from the exactness of the sequence (\ref{eq:coh0}).
\end{proof}

Let $A$ be a finite element system on $\calT$. Define a finite element system $\mathring{A}$ by:
\begin{equation}
\mathring{A}^k(T) = \{ u \in A^k(T) \ : \ \forall T'\in \calT \quad T' \subseteq T \ \Rightarrow \ u|_{T'} \textrm{ is } A\textrm{--harmonic} \}.
\end{equation}
We say that $\mathring{A}$ is the subsystem of locally harmonic forms.

\begin{proposition} \label{prop:harmiscomp} If $A$ is a compatible FE system then $\mathring{A}$ is a compatible FE system such that the De Rham map $\rho^k : \mathring{A}^k(\calT) \to \calC^k(\calT)$ is an isomorphism. 
\end{proposition}
\begin{proof}
This was essentially proved in \cite{Chr08M3AS}.
\end{proof}

This construction generalizes \cite{KuzRep05}, in which $\div$-conforming finite elements are defined on polyhedra in $\bbR^3$. Since $\calC^k(\calT)$ has a canonical basis, the De Rham map determines a corresponding canonical basis of $\mathring{A}^k(\calT)$. Its elements can be constructed by recursive harmonic extension.

\begin{example} On simplexes, lowest order trimmed polynomial differential forms  are locally harmonic in the sense of (\ref{eq:l2harmonic}), with respect to the $\rmL^2$ product associated with \emph{any} piecewise constant Euclidean metric \cite{Chr08M3AS}.
\end{example}

\begin{example} Locally harmonic forms can be used  to define finite element spaces on the dual cellular complex of a given simplicial one \cite{BufChr07}.  Choose a simplicial refinement of the dual mesh, for instance the barycentric refinement of the primal mesh. Consider Whitney forms on this refinement as a compatible finite element system on the dual mesh. Then take the subsystem of locally harmonic forms. In space dimension $d$, this provides a space of $k$-forms with the same dimension as the space of $(d-k)$-forms on the primal mesh. Duality in the sense of an inf-sup condition was proved in \cite{BufChr07} with applications to the preconditioning of integral operators appearing in electromagnetics. 
\end{example}

Duality methods are quite common in finite volume settings, for a recent development see \cite{AndBenKar10}.

\begin{example}
For a given fine mesh one can agglomerate elements into a coarser cellular mesh. The finite element system on the fine mesh then provides a finite element system on the coarse one, with identical function spaces. If the former system is compatible so is the latter. Associated with the the latter one can consider the subsystem of locally harmonic forms. This procedure can be applied recursively: at each level one can consider the locally harmonic forms of the finer level. This yields a multilevel analysis which can be used for multigrid preconditioning \cite{PasVas08}. 
\end{example}

Here is a third application, again involving the dual finite elements.
\begin{example}\label{ex:euler}
Recall that the incompressible Euler equation can be written:
\begin{equation}
\dot{u} + \div (u \otimes u) + \grad p = 0 \rmins{and} \div u = 0.
\end{equation}
Here, $u$ is a vector field with time derivative $\dot{u}$. One uses $\div$-conforming Raviart-Thomas elements for $u$ -- let $X_h$ denote this space. One uses a weak formulation with locally harmonic $\curl$-conforming elements on the dual grid as test functions -- let $Y_h$ denote this space. The $\rmL^2$ duality and its extensions to Sobolev spaces, is denoted $\langle \cdot , \cdot \rangle$. It is invertible on $X_h \times Y_h$, at least in dimension 2, as proved in \cite{BufChr07}. The semi-discrete problem reads as follows. Find a time dependent $u_h \in X_h$ such that, for all $v_h \in Y_h$ which are orthogonal to the subspace of $Y_n$ of $\curl$-free elements, we have: 
\begin{equation}
\langle \dot{u_h}, v_h \rangle + \langle \div (u_h \otimes u_h), v_h\rangle = 0.
\end{equation}
To see that the second bracket is well-defined, remark that $\div (u_h \otimes u_h)$ and $\div \div (u_h \otimes u_h)$ both have $\rmW^{-1,q}$ Sobolev regularity for all $q < 2$. On the other hand $v_h$ and $\curl v_h$ have $\rmW^{0,q'}$ regularity, for $q' >2$. 

Analogous spaces, with regularity expressed by a second order operator, have been used for Regge Calculus \cite{Chr11NM} and elasticity \cite{SchSin07}.
\end{example}

In computational fluid dynamics it is often important to include some form of upwinding in the numerical method to obtain stability. A model problem for this situation is the equation:
\begin{equation}\label{eq:convdiff}
\epsilon \Delta u + V \cdot \grad u = f,
\end{equation}
where the field $u \in \rmH^1(S)$ satifies homogeneous Dirichlet boundary condition, $f$ is a given forcing term, whereas $V$ is a given flow field, which we take to be divergence free. One is interested in the asymptotic behavior for small positive $\epsilon$ (viscosity). Such problems can for instance be solved with Petrov-Galerkin methods. Consider a cellular complex and a large compatible finite element system on it, obtained for instance by refinement. We let one space (say the trial space) consist of locally harmonic forms for the standard $\rmL^2$ product. For the other space (say the test space) we use locally harmonic forms for a weighted $\rmL^2$ product. 

To motivate our choice of weight we introduce some more notions of differential geometry. Given a $1$-form $\alpha$ we define the covariant exterior derivative:
\begin{equation}
\rmd_\alpha : u \mapsto \rmd u + \alpha \wedge u,
\end{equation}
These operators do not form a complex but we have:
\begin{equation}
\rmd_\alpha \rmd_\alpha u = (\rmd \alpha) \wedge u.
\end{equation}
In gauge theory the term $\rmd \alpha$ is called curvature. Supposing that $\alpha = \rmd \beta$ for a function $\beta$ we have:
\begin{equation}
\rmd_\alpha = \exp(-\beta) \rmd \exp (\beta) u.
\end{equation}
One says that $u \mapsto \exp (\beta) u$ is a gauge transformation.

We suppose the domain $S$ is equipped with a Riemannian metric. It provides in particular a $\rmL^2$ scalar product on differential forms. We denote by $\rmd_\alpha^\star$ the formal adjoint of $\rmd$. When $\alpha = \exp(-\beta)$ we have:
\begin{equation}
\rmd_\alpha^\star  u = \exp(\beta) \rmd^\star \exp(-\beta).
\end{equation} 
A natural generalization of (\ref{eq:l2harmonic}) is:
\begin{equation}
\rmd_\alpha^\star \rmd u = 0 \rmins{and} \rmd_\alpha^\star u = 0,
\end{equation}
which can be written:
\begin{equation}\label{eq:betaharm}
\rmd^\star \exp(-\beta) \rmd u = 0 \rmins{and} \rmd^\star \exp(-\beta) u = 0.
\end{equation}

\begin{example} To address (\ref{eq:convdiff}) define the $1$-form $\alpha$ by:
\begin{equation}
\alpha (\xi) = -\epsilon^{-1} V \cdot \xi.
\end{equation}
and remark that (\ref{eq:convdiff}) can be rewritten:
\begin{equation}
\rmd_\alpha^\star \rmd u = -\epsilon^{-1}f.
\end{equation}
Given a cellular complex $\calT$ on $S$, with flat cells, and a large compatible FE system $A$, we construct two spaces of locally harmonic forms, distinguished by the choice of scalar product $a$. For one (the trial space) take $a$ to be the $\rmL^2$ scalar product. For the other (the test space), we choose for each $T$, a constant approximation $\alpha_T$ of the pullback of $\alpha$ to $T$. Let $\beta_T$ be the affine function with zero mean on $T$ such that $\rmd \beta_T = \alpha_T$. For the trial spaces, use the scalar product defined on a cell $T$ by:
\begin{equation}\label{eq:weightscalp}
a(u,v) = \int_T \exp(-\beta_T) u \cdot v,
\end{equation}
to define the locally harmonic functions. If $v$ is a constant differential form, $u= \exp(\beta_T) v$ satisfies the equations (\ref{eq:betaharm}).

The canonical basis of the testspace will then be upwinded or downwinded (depending on the sign in front of $\beta_T$ in (\ref{eq:weightscalp})), compared with the canonical basis of the trial space.
\end{example}

\begin{example} Similar notions can be used to address the Helmholtz equation:
\begin{equation}\label{eq:helmholtz}
\Delta u + k^2 u = 0.
\end{equation}
We wish to construct a compatible FE system over $\bbC$, which contains a certain number of plane waves:
\begin{equation}
u_\xi : x \mapsto \exp(i \xi \cdot x).
\end{equation}
To  contain just one of them we remark that for any (flat) cell $T$, if $\xi_T$ is the tangent component of $\xi$ on $T$ we have, on $T$:
\begin{equation}
\Delta_T u_\xi|_T - i \xi_T \cdot \grad_T u_\xi|_T = 0.
\end{equation}
Connecting with the previous example, on a cell $T$ we let $\beta_T$ be the affine real function with zero mean and gradient $\xi_T$ on $T$. Define:
\begin{equation}
a(u,v) = \int_T \exp(-i\beta_T) u \cdot v.
\end{equation}
Here, extra care must be taken because this bilinear form is not positive definite and will in fact be degenerate at interior resonances of the cell. Away from them, the locally harmonic forms for the infinite dimensional element system $\hat X$ (with $q= 2$) are well behaved, in the sense of satisfying Proposition \ref{prop:harmiscomp}, and contain the plane wave $u_\xi$.

More generally, suppose one wants to construct a FE system containing a good approximation to a particular solution of (\ref{eq:helmholtz}). To the extent that the solution can be locally approximated by a plane wave, the finite element system will contain a good approximation of it,  for a choice of a family of functions $\beta_T$, one for each cell $T$, to be determined (maybe adaptively).
\end{example}


\section{Interpolators \label{sec:interpol}}

\paragraph{Mirrors and interpolators.} The notion of mirror system formalizes that of degrees of freedom, with particular emphasis on their geometric location.
\begin{definition} 
A \emph{mirror system} is a choice, for each $k$ and $T$, of a subspace $\calZ^k(T)$ of $\hat{X}^k(T)^\star$, called a $k$-mirror on $T$.
\end{definition} 
Any $k$-form $u$ in $\hat X^k(T)$ then gives a linear form $\langle \cdot , u \rangle$ on $\calZ^k(T)$ which we call the mirror image of $u$. For a global $k$-form $u$ the mirror images can be collected into a single object. We define:
\begin{equation}
\Phi^k u = \langle \cdot, u|_T \rangle_{T \in \calT} \in \calZ^k(\calT)^\star,
\end{equation}
where, for any subcomplex $\calT'$ of $\calT$,  we define a (global) mirror $\calZ^k(\calT')$ by:
\begin{equation}
 \calZ^k(\calT')=\bigoplus_{T \in \calT'} \calZ^k(T).
\end{equation}
We say that a mirror system is \emph{faithful} to an element system $A$ if, for any subcomplex $\calT'$, restricting $\Phi$ to $\calT'$ determines an isomorphism:
\begin{equation}
\Phi^k(\calT') : A^k(\calT') \to \calZ^k(\calT')^\star.
\end{equation}

\begin{example}
The canonical mirror system for the trimmed polynomial FE system of order $p$ is the following, where $\dim T = m$.
\begin{equation}
\calZ^k(T) = \{ u \mapsto \ts \int_T v \wedge u  \  :  \  v \in \poly\alt_{p - m + k -1}^{m -k}(T) \}.
\end{equation}
It follows from results in \cite{ArnFalWin06} that it is faithful.
\end{example}

\begin{proposition}\label{prop:faithfuliff}
When $A$ admits extensions, a given mirror system $\calZ$ is faithful iff the duality product on $ \calZ^k(T) \times A^k_0(T) $ is invertible for each $k$ and $T$.
\end{proposition}
\begin{proof}
(i) Suppose $\calZ$ is faithful. By induction on dimension we get $\dim \calZ^k(T) = \dim A^k_0(T)$. Moreover $\Phi^k(T)$ induces an injection $A^k_0(T) \to \calZ^k(T)^\star$. Thus duality on $ \calZ^k(T) \times A^k_0(T) $ is invertible.

(ii) Suppose duality on $ \calZ^k(T) \times A^k_0(T) $ is invertible for all $k$ and $T$. Then:
\begin{equation}
\dim A^k(\calT') = \ts \sum_{T \in \calT'} \dim A^k_0(T) = \dim \calZ^k(\calT')^\star.
\end{equation}
Moreover $\Phi^k(\calT')$ is injective. Indeed if $\Phi^k u = 0$ then for $T \in \calT$, $u|_T$ is proved to be $0$ by starting with cells $T$ of dimension $k$, and incrementing cell dimension inductively using $u|_{\partial T} = 0$.  
\end{proof}

\begin{definition}
For a finite element system $A$, an \emph{interpolator} is a collection of projection operators $I^k(T) : \hat X^k(T) \to A^k(T)$, one for each $k \in \bbN$ and $T \in \calT$, which commute with restrictions to subcells.
\end{definition}
One can then denote it simply with $I^\bs$ and extend it unambiguously to any subcomplex $\calT'$ of $\calT$. Any faithful mirror system defines an interpolator by:
\begin{equation}
\Phi^k I^k u = \Phi^k u.
\end{equation}
We call this the interpolator associated with the mirror system.

\begin{proposition}The following are equivalent:
\begin{itemize}
\item $A$ admits extensions,
\item $A$ has a faithful mirror system,
\item $A$ can be equipped with an interpolator.
\end{itemize}
\end{proposition}
\begin{proof}
(i) Suppose $A$ admits extensions. For each $k$ and $T$, choose a closed supplementary of $A^k_0(T)$ in $\hat X^k(T)$ and let $\calZ^k(T)$ be its annihilator. The duality product on $ \calZ^k(T) \times A^k_0(T) $ is then invertible for each $k$ and $T$, so that $\calZ$ is faithful to $A$.

(ii) As already indicated, any faithful mirror system defines an interpolator.

(iii) Suppose $A$ has an interpolator. If $u \in A^k(\partial T)$ extend it in $\hat X^k(T)$ and interpolate it to get an extension in $A^k(T)$.
\end{proof}

\begin{example} In particular the trimmed polynomial FE system of order $p$ is compatible. As remarked in \cite{Chr10CRAS}, it is minimal among compatible finite element systems containing polynomial differential forms of order $p-1$.
\end{example}

For a given element system $A$, admitting extensions, it will be useful to construct extension operators $A^k(\partial T) \to A^k(T)$, i.e. linear left inverses of the restriction. One also remarks that a faithful mirror system determines a particular extension. Namely, to $u\in A^k(\partial T)$ one associates the unique $v \in A^k(T)$ extending $u$ and such that for all $l\in \calZ^k(T)$,  $l(v) = 0$.

\begin{proposition} Let $A$ be an element system on $\calT$ admitting extensions. Let $\calM$ denote the set of mirror systems that are faithful to $A$, $\calI$ the set of  interpolators onto $A$ and $\calE$ the set of extensions in $A$. The natural map $\calM \to \calI \times \calE$ is bijective.
\end{proposition}
\begin{proof}
For a given interpolator $I\in \calI$ and $E\in \calM$, define a mirror system $z(I,E)$ as follows. For $k\in \bbN$ and $T \in \calT$ consider the map:
\begin{equation}
Q^k_T = (\id - E \partial) \circ I \ : \ \hat{X}^k(T) \to \hat{X}^k(T). 
\end{equation}
It is a projector with range $A^k_0(T)$. Let $z(I,E)^k(T)$ denote the annihilator of its kernel. In other words $l \in z(I,E)^k(T)$ iff:
\begin{equation}
\forall u \in \hat{X}^k(T) \quad I u = E I \partial u  \ \Rightarrow \ l(u) = 0.
\end{equation}
We claim that $z$ inverts the given map $a : \calM \to \calI \times \calE$. 

(i) Given $(I,E)$ we check that the interpolator and extension defined by $z(I,E)$ are $I$ and $E$.

-- Pick $u \in \hat{X}^k(T)$.  We have $I(u - Iu)= 0$ and $E I \partial (u - Iu) =0$ hence $l(u- Iu) = 0$. The interpolator deduced from $z(I,E)$ is thus $I$.

-- Pick $u \in A^k(\partial T)$. We have $ I E u = E I \partial E u$ hence for all $l\in z(I,E)^k(T)$ we have $l(u)= 0$. The extension deduced from $z(I,E)$ is thus $E$.

(ii) Given a faithful mirror system $\calZ$ with associated interpolators $I$ and $E$, we check that $z(I,E) = \calZ$.

Pick $l \in \calZ^k(T)$. If $u \in \hat{X}^k(T)$ is such that $I u = E I \partial u$ then $l(u) = l(Iu) = l(EI \partial u) = 0$. Hence $l \in  z(I,E)^k(T)$. On the other hand $z(I,E)$ is also faithful to $A$ by Proposition \ref{prop:faithfuliff}. The inclusion $\calZ^k(T)  \subseteq z(I,E)^k(T)$ then implies equality.
\end{proof}

\begin{example} Let  $A$ be a finite element system. Equip each $\hat X^k(T)$ with a continuous bilinear form $a$ which is non-degenerate on $A^k_0(T)$ (e.g. $\hat X^k(T)$ is continuously embedded in a Hilbert space).  One can define a mirror system by:
\begin{equation}
\calZ^k(T)= \{a( \cdot, v) \ : v \in A^{k}_0(T)\}.
\end{equation}
When the FES admits extensions, the associated interpolator can be interpreted as a recursive $a$-projection, starting from cells of minimal dimension, and continuing by incrementing dimension at each step, projecting with respect to $a$ with given boundary conditions.
\end{example}

\paragraph{Commuting interpolators. }
It is of interest to construct interpolators which commute with the exterior derivative. When $T' \subseteq T$, restriction maps $\hat X^k(T) \to \hat X^k(T')$, so that a $k$-mirror $\calZ^k(T')$ on $T'$ can also be considered as a $k$-mirror on $T$. If the mirror system is faithful it must be in direct sum with  $\calZ^k(T)$. 
\begin{proposition}
An interpolator commutes with the exterior derivative if and only if its mirror system satisfies:
\begin{equation}\label{eq:dadjoint}
\forall l \in \calZ^k(T) \quad  l \circ \rmd \in \calZ^{k-1}(\tilde{T}).
\end{equation}
\end{proposition}
\begin{proof}
(If) Pick $u\in \hat X^{k-1}(T)$.  Suppose first $\Phi^{k-1} u = 0$. Then for all $l \in \calZ^k(T)$, $l (\rmd u) = 0$, hence $\Phi^k(\rmd u) = 0$. In the general case, since $\Phi^{k-1}(u - \Phi^{k-1} u) = 0$ we deduce $\Phi^k(\rmd u - \rmd \Phi^{k-1}u) = 0$, so that $\Phi^k \rmd u = \rmd \Phi^{k-1}u$.

(Only if) Suppose $l \in \calZ^k(T)$ and $ l \circ \rmd \not \in \calZ^{k-1}(\tilde{T})$. Pick $u\in \hat X^{k-1}(T)$ such that $l(\rmd u) \neq 0$ but for all $l'\in \calZ^{k-1}(\tilde{T})$, $l'(u)=0$. Then $I^k \rmd u \neq 0$ but $I^{k-1} u = 0$ (so $\rmd I^{k-1} u = 0$).
\end{proof}

\begin{example}
The canonical mirror system of trimmed polynomials of order $p$ yields a commuting interpolator.
\end{example}

Suppose that $\calZ$ is a mirror system on $\calT$ such that (\ref{eq:dadjoint}) holds. Then we have a well defined map $\hat \rmd: l\mapsto l \circ \rmd $ from $\calZ^k(\calT)$ to $\calZ^{k-1}(\calT)$. Denote by $\delta$ its adjoint, which maps from $\calZ^{k-1}(\calT)^\star$ to $\calZ^k(\calT)^\star$.
\begin{remark} The following diagram commutes:
\begin{equation}
\xymatrix{
\hat X^{k-1}(\calT) \ar[r]^\rmd \ar[d]^{\Phi^{k-1}} & \hat X^{k}(\calT) \ar[d]^{\Phi^k}\\
\calZ^{k-1}(\calT)^\star \ar[r]^\delta & \calZ^k(\calT)^\star
}
\end{equation}
\end{remark}
\begin{proof}
Pick $l \in \calZ^k(\calT)$ and $u \in \hat X^{k-1}(\calT)$. We have:
\begin{equation}
(\delta \Phi u)(l)= (\Phi u)(\hat \rmd l) =(\Phi u) (l\circ \rmd) = l(\rmd u)= \Phi(\rmd u).
\end{equation}
This concludes the proof.
\end{proof}

\begin{proposition}
Equip each $\hat X^k(T)$ with a continuous scalar product $a$. For any compatible finite element system $A$, the following is a faithful mirror system yielding a commuting interpolator:
For $k = \dim T$:
\begin{equation}
\calZ^k(T) = \{a( \cdot, v) \ : v \in \rmd A^{k-1}_0(T)\} + \{ \bbR \ts \int \cdot \}
\end{equation}
For $k < \dim T$:
\begin{equation}
\calZ^k(T) = \{a( \cdot, v) \ : v \in \rmd A^{k-1}_0(T)\} + \{ a( \rmd \cdot, v) \ : v \in \rmd A^{k}_0(T)\}.
\end{equation}
 \end{proposition}
This is the natural generalization, to the adopted setting, of projection based interpolation, as defined in \cite{DemBab03}\cite{DemBuf05}. When the scalar products $a$ are all the $\rmL^2$ product on forms, we call it \emph{harmonic interpolation}.

Suppose $\calU$ is a cellular complex on $M$ and $\calV$ is a cellular complex on $N$  giving rise to the product complex $\calU \times \calV$ on $M \times N$. If $\calZ$ is a mirror system on $\calU$ and $\calY$ is a mirror system on $\calV$, then $\calZ \otimes \calY$ is a mirror system on $\calU \times \calV$ defined by:
\begin{equation}
(\calZ \otimes \calY)^k(U \times V) = \bigoplus_l \calZ^l(U) \otimes \calY(V)^{k-l}.
\end{equation}

\begin{proposition} The tensor product of two faithful mirror systems is faithful.
\end{proposition}
\begin{proof}
We have three FE systems: $A$, $B$ and $A \otimes B$, all admitting extensions.
Moreover by Proposition \ref{prop:czero}:
\begin{equation}
(A \otimes B)^k_0(U \times V) = \bigoplus_l A^l_0(U) \otimes B^{k-l}_0(V).
\end{equation}
Now apply Proposition \ref{prop:faithfuliff}.
\end{proof}

\paragraph{Extension-projection interpolators.}
It will also be of interest to construct extension operators which commute with the exterior derivative. More precisely for a cell of dimension $m$ a commuting extension operator is a family of  operators $E^k: A^k(\partial T) \to A^k(T)$ for $0\leq k\leq m-1$, together with a map $E^m : \bbR \to A^m(T)$ such that the following diagram commutes:
\begin{equation}\label{eq:comext}
\xymatrix{
0 \ar[r] & \bbR \ar[r]    \ar[d]         & A^0 (\partial T)    \ar[r]  \ar[d]^{E^0}& \cdots \ar[r]    \ar[d]     & A^{m-1}(\partial T) \ar[r] \ar[d]^{E^{m-1}} & \bbR               \ar[d]^{E^m}\ar[r] &  0\\
0 \ar[r] & \bbR \ar[r]\  & A^0(T)          \ar[r] & \cdots \ar[r] & A^{m-1}(T)     \ar[r]   & A^m(T) \ar[r] & 0
}
\end{equation}
When $A^{m-1}(\partial T)$ has an element with non-zero integral an $E^m$ such that the above diagram commutes, is uniquely determined by $E^{m-1}$ and exists iff $E^{m-1}$ sends elements with zero integral to closed forms.

\begin{proposition}\label{prop:extcom} Let $A$  be a finite element system where each  $A^k(T)$ is equipped with a scalar product $a$. Suppose $T$ is a cell such that (\ref{eq:coh0}) is exact. Put $m =\dim T$. If $T$ admits extensions, the harmonic extension operators defined by Proposition \ref{prop:intext} commute in the sense of diagram (\ref{eq:comext}).
\end{proposition}

We suppose that for each cell $T$, extension operators $E : A^\bs(\partial T) \to A^\bs(T)$  and projections  $P: X^\bs(T) \to A^\bs(T)$ have been defined. Then $E$ and $P$ uniquely determine an interpolator $J$ as follows. One constructs $J$ inductively. The initialization on cells of dimension $0$ (points) is trivial. Let now $T$ be a cell and suppose that $J$ has been constructed for all cells on its boundary. On $T$ we define, for $u\in \hat X^k(T)$ with $k < \dim T$:
\begin{align}
J u & = E J \partial u   + (\id - E \partial) Pu,\\
    & = Pu + E(J \partial u -\partial Pu).
\end{align}
and for $k = \dim T$ we simply put:
\begin{equation}
J u = Pu. 
\end{equation}
The only thing to check is that $J$ commutes with restriction from $T$ to cells on the boundary, which is trivial. We call $J$ the associated extension-projection (EP) interpolator.

Let $T$ be cell of dimension $m$. We say that an endomorphism $F$ of $X^m(T)$ preserves integrals if $\int F u = \int u$ for all $u \in X^m(T)$. An interpolator is said to preserve integrals if it preserves integrals on all cells. For an interpolator $I$, this is equivalent to commutation of the following diagram, involving De Rham maps:
\begin{equation}
\xymatrix{
\hat X^\bs(\calT) \ar[rr]^I \ar[dr]_\rho& & A^\bs(\calT) \ar[dl]^\rho \\
 & \calC^\bs(\calT)
}
\end{equation}

We would like furthermore the interpolator to commute with the exterior derivative.
\begin{proposition}
Suppose that the projectors $P: X^\bs(T) \to A^\bs(T)$ commute with the exterior derivative, that the extensions $E$ commute in the sense of diagram (\ref{eq:comext}) and that $P$ preserves integrals. Then the associated EP interpolator commutes with the exterior derivative and preserves integrals.
\end{proposition}
\begin{proof}
By induction on the dimension of cells. Let $T$ be a cell of dimension $m$.

If $u \in \hat X^k(T)$ with $k\leq m-2$ the commutation $\rmd J u = J \rmd u$ follows immediately from the commutation of $P$,  $E$ and $\partial$,  as well as $J$ on $\partial T$ (which follows from the induction hypothesis). 

For $u \in \hat X^{m-1}(T)$we have:
\begin{equation}
 \rmd J u = \rmd P u + \rmd E ( J \partial u- \partial Pu).
\end{equation}
But we have:
\begin{align}
\int_{\partial T}  ( J \partial u- \partial Pu) & = \int_{\partial T} \partial u - \int_T \rmd Pu  = \int_{\partial T} \partial u  - \int_T P \rmd u,\\
 & = \int_{\partial T} \partial u  - \int_T  \rmd u= 0.
\end{align}
The first thing we used is that $J$ preserves integrals on the boundary (which is an induction hypothesis). Hence by (\ref{eq:comext}):
\begin{equation}
\rmd E ( J \partial u- \partial Pu) = 0.
\end{equation}
Thus:
\begin{equation}
\rmd J u = \rmd P u = P \rmd u = J \rmd u.
\end{equation}

That $J$ preserves integrals on $T$ follows simply from the fact that, for $u \in X^m(T)$, we have $\int J u = \int P u = \int u$.
\end{proof}

Compared with the use of mirror systems, the advantage of defining an interpolator from extensions and projections is that approximation properties of the interpolator follow directly from estimates on the extensions and projections. In the following  we denote the $\rmL^q(U)$ norm simply by $\| \cdot \|_U$.
\begin{proposition} \label{prop:orderopt} Let $(\calT_n)$ be a sequence of cellular complexes, each equipped with a compatible FE system. Suppose that for each $T$, $\tilde X^k(T)$ is equipped with a densely defined seminorm $\llp \cdot \rrp_T$. We also suppose that we have functions  $\lambda_n:\calT_n \to \bbR^\ast_+$ and $\tau_n: \calT_n \to \bbR^\ast_+$ and extensions and projections satisfying:
\begin{align}
\|u -P_n u\|_{T} & \cleq \tau_n(T) \llp  u \rrp_{T},\\
\|u -P_n u\|_{\partial T} & \cleq \tau_n(T) \lambda_n(T)^{-1} \llp u\rrp_{T},\\
\|E_n u\|_{T}  & \cleq  \lambda_n(T) \|u \|_{\partial T}.
\end{align}
We suppose that for $T,T' \in \calT_n$, if $T' \subseteq T$ then  $\lambda_n(T')\ceq \lambda_n(T)$ and $\tau_n(T')\ceq \tau_n(T)$.
Then the associated EP interpolator satisfies:
\begin{equation}
\|u  - I_n u\|_{T} \cleq  \tau_n(T) \sum_{T' \subseteq T}\lambda_n(T)^{\dim T - \dim T'}  \llp u \rrp_{T'},
\end{equation}
where we sum over subcells $T'$ of $T$ in $\calT_n$.
\end{proposition}

We shall express this bound as order-optimality.
\begin{proof} By induction.
Suppose $T$ is a cell and that $I_n$ is order-optimal on its boundary. For $u\in \hat X^k(T)$ with $k < \dim T$, we then have:
\begin{align}
\|u -I_nu\|_T & \cleq \|u-P_nu\|_T + \lambda_n(T) \| I_n\partial u-\partial P_nu\|_{\partial T} ,\\
& \cleq \|u-P_nu\|_T + \lambda_n{(T)} \| u-P_nu\|_{\partial T} + \lambda_n(T) \| u- I_nu\|_{\partial T},\\
& \cleq \tau_n(T) \llp u \rrp + \lambda_n(T) \| u- I_nu\|_{\partial T}.
\end{align}
This completes the proof.
\end{proof} 

This proposition was designed with the $p$--version of the finite element method in mind. One can think of $\tau_n(T)= p^{-1}$ and $\lambda_n(T)= p^{-1/q}$. The seminorms involved would correspond to Sobolev spaces, possibly weighted. Thus one would require extension operators whose $\rmL^q(\partial T) \to\rmL^q(T)$ norm is of order $p^{-1/q}$. The construction used to prove Proposition 3.3 in the preprint of \cite{Chr07NM} might be useful here. 



\section{Quasi-interpolators \label{sec:quasiint}}

Let $S$ be a domain in $\bbR^d$. Let $(\calT_n)$ be a sequence of cellular complexes on $S$, each equipped with a compatible FE system $A[n]$. We define $X^k_n= A[n]^k(\calT_n)$. 

\begin{definition} Consider a sequence of maps $Q^k_n: \rmL^q(S) \to X^k_n$. 
\begin{itemize}
\item We say that they are \emph{stable} if for $u \in \rmL^q(S)$:
\begin{equation}
\| Q^k_nu  \|_{\rmL^q(S)} \cleq \|u\|_{\rmL^q(S)} ,
\end{equation}

\item We say that they are \emph{order optimal} if for $u\in \rmW^{\ell,q}(S)$:
\begin{align}
\| u - Q^k_nu \|_{\rmL^q(S)}  \cleq \tau_n^\ell \| u\|_{\rmW^{\ell,q}(S)}. 
\end{align}
where $\tau_n^\ell$ is the order of best approximation on $X^k_n$ in $\rmW^{\ell,q}(S) \to \rmL^q(S)$ norm.

\item We say that they are \emph{quasi-projections} if for some $\alpha < 1$ we have for all $u \in X^k_n$:
\begin{equation}
\| u - Q^k_nu \|_{\rmL^q(S)} \leq \alpha \| u\|_{\rmL^q(S)}.
\end{equation}

\item We say that they \emph{commute} if they commute with the exterior derivative. 
\end{itemize}
\end{definition}

We shall construct stable quasi-interpolators of the form:
\begin{equation}
Q= I R E.
\end{equation}
where $I$ is an interpolator, $R$ is a regularization (smoother) approximating the identity and $E$ is an extension operator (usually $R$ requires values outside $S$).

The following technique will be referred to as \emph{scaling}. 
\begin{lemma}\label{lem:scaling}
For a cell $T$ of diameter $h_T$ and barycentre $b_T$, consider the scaling map $\sigma_T: x \to h_T x + b_T$ and let $\hat T$ be the preimage of $T$ by $\sigma_T$, called the reference cell.
Let $m$ be the dimension of $T$ ($0 \leq m \leq d$). Let $u$ be a $k$-multilinear form on $T$ and $\hat u = \sigma_T^\star u$ the pull-back of $u$ to $\hat T$. Then we have:
\begin{equation}
\| u \|_{\rmL^q(T)} = h_T^{-k + m/q} \|\hat u \|_{\rmL^q(\hat T)}.
\end{equation} 
\end{lemma}

We adopt the $h$-setting, in that the differential elements $A[n]^k(T)$, when pulled back to reference domains $\hat T$, belong to compact families, see \cite{ArnFalWin06} Remark p. 64. This hypothesis excludes for instance methods where the polynomial degree is unbounded as $n \to \infty$ and normally requires the cells to be shape-regular.

\begin{proposition} \label{prop:orderoptimal}
Suppose the finite elements $A[n]^k(T)$ contain polynomials of degree $\ell - 1$. A sequence of commuting interpolators $I_n$ can be constructed to satisfy:
\begin{equation}\label{eq:orderopt}
\|u  - I_n u\|_{\rmL^q(T)} \cleq  \sum_{T' \subseteq T} h_T^{\ell+(\dim T - \dim T')/q}  \| \nabla^\ell u\|_{\rmL^q(T')},
\end{equation}
where we sum over subcells $T'$ of $T$ in $\calT_n$.
\end{proposition}
We shall express this bound as order-optimality.
\begin{proof}
Choose commuting interpolators $I_n$, which when pulled back to the reference cell satisfy:
\begin{equation}
\|u  - \hat I_n u\|_{\rmL^q(\hat T)} \cleq  \sum_{T' \subseteq T}  \| \nabla^\ell u\|_{\rmL^q(\hat T')}.
\end{equation}
Estimate (\ref{eq:orderopt}) follows by scaling. Such interpolators can be constructed from mirror systems (e.g. harmonic interpolation), or from extensions and projections as in Proposition \ref{prop:orderopt}.
\end{proof} 

\paragraph{Regularizer}

We consider regularizing operators constructed as follows. We require a function $\psi$ supported in the unit ball $\bbB^d$ and a function:
\begin{equation}
\Phi: \left \{\begin{array}{l}
\bbB^d \times \bbR^d \to \bbR^d,\\
 (y,x) \mapsto \Phi_y(x).
\end{array}
\right.
\end{equation}
For any given $y$ we denote by $\Phi_y^\star$ the pullback by $\Phi_y : \bbR^d \to \bbR^d$. That is for any $k$-multilinear form $u$:
\begin{equation}
(\Phi_y^\star u)[x](\xi_1, \ldots , \xi_k)= u[\Phi_y(x)](\rmD_x \Phi_y(x) \xi_1, \ldots, \rmD_x \Phi_y(x)\xi_k).
\end{equation}
We define:
\begin{equation}\label{eq:defreg}
R u = \int_{\bbB^d} \psi(y) \Phi_y^\star u \, \rmd y
\end{equation}
For the function $\psi$ we choose one which is smooth, rotationally invariant, non-negative, with support in the unit ball and with integral 1. Later on, we will require convolution by $\psi$ to preserve polynomials up to a certain degree, see (\ref{eq:polypres}), but until then this hypothesis will be irrelevant. In what follows, when we integrate with respect to the variable $y$ it will always be on the unit ball $\bbB^d$ so that we omit this from expressions such as (\ref{eq:defreg}).

Given such a $\psi$ we are interested in how properties of $\Phi$ reflect upon $R$. For this purpose we call $R$ defined by (\ref{eq:defreg}) the regularizer associated with $\Phi$. To emphasize its dependence on $\Phi$ we sometimes denote the regularizer by $R[\Phi]$. Formally, $R[\Phi]$ will commute with the exterior derivative, since pullbacks do. The minimal regularity we assume of $\Phi$ is to have continuous second derivatives. We also suppose that for given $x$, $y \mapsto \Phi_y(x)$ is a diffeomorphism from $\bbB^d$  to its range. This is enough for the commutation to hold. In what follows we shall be a bit more careful about the regularizing effects of $R[\Phi]$ for given $\Phi$. 

Let $\bbB(x,\delta)$ denote the ball with center $x$ and radius $\delta$. For any subset $U$ of $\bbR^d$ its $\delta$-neighborhood is defined as:
\begin{equation}
\calV^\delta(U) = \cup \{\bbB(x,\delta) \ : \ x \in U\}.
\end{equation}

As a first result we state:
\begin{proposition}\label{prop:regbound}
For any $\delta >0 $ and $C > 0$, there exists $C' > 0$ such that, for all $\Phi$ satisfying, at a point $x$, for all $y \in \bbB^d$:
\begin{equation}
| \Phi_y (x) -x| \leq \delta,
\end{equation}
and: 
\begin{equation}
\| \rmD_y \Phi_y(x)^{-1} \|,\, \| \rmD_x \Phi_y(x)\| \leq C,
\end{equation}
the associated regularizer satisfies:
\begin{equation}\label{eq:regbound}
\|(R[\Phi] u)(x)\| \leq C' \| u\|_{\rmL^1(\calV^\delta(x))}.
\end{equation}
\end{proposition}
\begin{proof}
With $x$ fixed, the Jacobian of the map:
\begin{equation}
\bbB^d \ni y \mapsto \Phi_y(x) \in \calV^\delta(x),
\end{equation}
has an inverse bounded by $C$.
\end{proof}

We are interested in estimates on derivatives of $R[\Phi] u$ when $u$ is a $k$-form. We have:
\begin{equation}\label{eq:nablaruexpr}
(\nabla R[\Phi] u)(x) = \int \psi(y) \rmD_x (\Phi_y^\star u)[x]  \, \rmd y,
\end{equation}
where we can substitute the expression:
\begin{align}\label{eq:dxphiu}
 &\rmD_x \Phi_y^\star u [x](\xi_0, \ldots, \xi_k) \nonumber\\
 =&  \rmD u [\Phi_y(x)](\rmD_x \Phi_y(x) \xi_0, \ldots, \rmD_x \Phi_y(x)\xi_k) \nonumber \\
 & + \ts \sum_{i=1}^k u [\Phi_y(x)](\rmD_x \Phi_y(x) \xi_1, \ldots, \rmD_{xx}^2 \Phi_y(x)(\xi_0, \xi_i), \ldots , \rmD_x \Phi_y(x)\xi_k).
\end{align}
The purpose of the following lemma is to get an integral expression for $\nabla R[\Phi] u$ not involving any derivatives of $u$. Essentially in (\ref{eq:nablaruexpr}), derivatives acting on $u$ are transferred to other terms under the integral sign, using integration by parts. Without the integral sign, this corresponds to identifying the derivatives of $u$ as a ``total'' divergence, up to expressions involving no derivatives of $u$.
\begin{lemma}
We have:
\begin{align}
&\rmD_x \Phi_y^\star u [x](\xi_0, \ldots, \xi_k) \nonumber \\
 = & \ts \sum_{i=1}^k u [\Phi_y(x)](\rmD_x \Phi_y(x) \xi_1, \ldots, \rmD_{xx}^2 \Phi_y(x)(\xi_0, \xi_i), \ldots , \rmD_x \Phi_y(x)\xi_k) \nonumber \\
 & -\ts  \sum_{i=1}^k u [\Phi_y(x)](\rmD_x \Phi_y(x) \xi_1, \ldots, \rmD_{yx}^2 \Phi_y(x)(\Xi_y(x)\xi_0, \xi_i), \ldots , \rmD_x \Phi_y(x)\xi_k) \nonumber \\
 & + \rmD_y \Phi_y^\star u [x](\Xi_y(x)\xi_0, \xi_1, \ldots, \xi_k),
\end{align}
with $\Xi_y(x)$ defined by:
\begin{equation}
\Xi_y(x) \xi = \rmD_y \Phi_y(x)^{-1} \rmD_x \Phi_y(x) \xi.
\end{equation}
We also have:
\begin{align}
  & \psi(y)\rmD_y \Phi_y^\star u [x](\Xi_y(x)\xi_0, \xi_1, \ldots, \xi_k) \nonumber \\
= & \div_y\big(\psi(y)\Xi_y(x)\xi_0 \, \Phi_y^\star u [x]( \xi_1, \ldots, \xi_k)\big) \nonumber\\
 & - \div_y\big(\psi(y)\Xi_y(x)^\star\big) \xi_0 \, \Phi_y^\star u [x]( \xi_1, \ldots, \xi_k).
\end{align}
\end{lemma}

The following proposition shows that the regularizer maps forms of $\rmL^1$ regularity to continuously differentiable forms, while being careful about the operator norm.
\begin{proposition}\label{prop:nablareg}
For any $\delta >0 $ and $C > 0$, there exists $C' > 0$ such that, for all $\Phi$ satisfying, at a point $x$, for all $y \in \bbB^d$:
\begin{equation}
| \Phi_y (x) -x| \leq \delta,
\end{equation}
and: 
\begin{equation}
\left. \begin{array}{c}
\| \rmD_y \Phi_y(x)^{-1} \|,\, \| \rmD_x \Phi_y(x)\|,\, \| \rmD_{xx}^2 \Phi_y(x) \| \\
\| \rmD_{yx}^2 \Phi_y(x) \|,\,  \| \rmD_y (\rmD_y \Phi_y(x)^{-1}) \|
\end{array} \right \} \leq C,
\end{equation}
the associated regularizer satisfies:
\begin{equation}\label{eq:nablalone}
\|(\nabla R[\Phi] u)(x)\| \leq C' \| u\|_{\rmL^1(\calV^\delta(x))}.
\end{equation}
\end{proposition}
\begin{proof}
The preceding Lemma gives an expression for $(\nabla R[\Phi] u)(x)$, from which the claim follows.
\end{proof}

In the estimate (\ref{eq:nablalone}) one has one order more of differentiation on the left hand side. In the following proposition we consider an equal amount of differentiation on both sides. 
\begin{proposition}\label{prop:nablastab}
Pick an integer $\ell \geq 0$, and  $\delta >0 $. For any $C > 0$, there exists $C' > 0$ such that, for all $\Phi$ satisfying, at a point $x$, for all $y \in \bbB^d$:
\begin{equation}
| \Phi_y (x) -x| \leq \delta,
\end{equation}
\begin{equation}\label{eq:jacest}
\| \rmD_y \Phi_y(x)^{-1} \| \leq C,
\end{equation}
and: 
\begin{equation}\label{eq:bder}
\| \rmD_x \Phi_y(x)\|,\, \| \rmD_{xx}^2 \Phi_y(x) \| , \ldots, \| \rmD^{\ell+1}_{x \ldots x}\Phi_y(x)\|  \leq C,
\end{equation}
the associated regularizer satisfies:
\begin{equation}
\|(\nabla^{\ell} R[\Phi] u)(x)\| \leq C' \| u\|_{\rmW^{\ell,1}(\calV^\delta(x))}.
\end{equation}
\end{proposition}
\begin{proof}
For $\ell = 0$ one uses the change of variable formula in the definition (\ref{eq:defreg}), the Jacobian being taken care of by estimate (\ref{eq:jacest}). The case $\ell = 1$ follows from expression (\ref{eq:dxphiu}). Differentiating this expression more times gives the claimed result for general $\ell$ (an expression for the differential  of order $\ell$ of a pullback will be given later).
\end{proof}

From now on we suppose $\Phi$ has the form:
\begin{equation}\label{eq:formphi}
\Phi_y(x) = x + \phi(x)y,
\end{equation}
for some function $\phi: S \to \bbR$. Most of the discussion would work with matrix valued maps $\phi : S \to \bbR^{d \times d}$ which could be important for anisotropic meshes. But for simplicity we do not consider anisotropic meshes here, and then scalar $\phi$ are sufficient. We also suppose that $\psi$ satisfies:
\begin{equation}\label{eq:polypres}
\int \psi(y) f(y) \rmd y = f(0),
\end{equation}
for all polynomials $f$ of degree at most $p+d$, for some integer $p \geq 0$.

The purpose of these hypotheses is to make the regularizer preserve polynomials of degree up to $p$. To see this, remark first that property (\ref{eq:polypres}) guarantees that convolution by $\psi$ preserves polynomials of degree $p+d$. We state:
\begin{proposition} 
Suppose $|\phi(x)| \leq \delta$ for some $\delta >0$.

If $u$ is also a polynomial of degree at most $p$ on  $\calV^\delta(x)$, then $\nabla^{\ell} Ru(x) = \nabla^\ell u(x)$ for all $\ell$.
\end{proposition}
\begin{proof}
We have:
\begin{equation}
\rmD_x \Phi_y(x) \xi= \xi + (\rmD \phi(x) \xi)y.
\end{equation}
Suppose $u$ is a $k$-form which is polynomial of degree at most $p$. Recall that:
\begin{equation}
(\Phi_y^\star u)[x](\xi_1, \ldots , \xi_k)= u[\Phi_y(x)](\rmD_x \Phi_y(x) \xi_1, \ldots, \rmD_x \Phi_y(x)\xi_k).
\end{equation}
As a function of $y$ this is a polynomial of degree at most $p + k \leq p + d$. Its value at $y=0$ is:
\begin{equation}
u[x](\xi_1, \ldots, \xi_k).
\end{equation}
This gives the case $\ell = 0$ of the proposition. For $\ell= 1$ one uses expression (\ref{eq:dxphiu}). Greater $\ell$ are obtained by further differentiation of this expression.
\end{proof}

To see to what extent Propositions \ref{prop:regbound}, \ref{prop:nablareg} and \ref{prop:nablastab} can be applied we remark: 
\begin{align}
\Phi_y(x) -x & = \phi(x)y,\\
\rmD_y \Phi_y(x) \xi & = \phi(x)\xi,\\
\rmD_y (\rmD_y \Phi_y(x)^{-1}) & = 0,\\
\rmD_x \Phi_y(x) \xi & = (\rmD \phi(x) \xi) y,\\
\rmD_{yx}^2 \Phi_y(x) (\xi, \xi') & = (\rmD \phi(x) \xi) \xi',\\
\rmD_{x \ldots x}^{\ell} \Phi_y(x)(\xi_1, \ldots, \xi_\ell) & = \rmD^\ell \phi(x)(\xi_1, \ldots, \xi_\ell)\, y.
\end{align}

\begin{proposition}\label{prop:nablastabopt}
Pick $\ell \leq p+1$. For any $\delta >0$ and any $C >0$  there exists $C' > 0$ such that, for all $\phi$ satisfying, at some point $x$:
\begin{equation}
|\phi(x)| \leq \delta,
\end{equation}
\begin{equation}
| \phi(x)^{-1} | \leq C,
\end{equation}
and:
\begin{equation}
\| \rmD_x \phi(x)\|,\ \| \rmD_{xx}^2 \phi(x) \| , \ldots, \| \rmD^{\ell+1}_{x \ldots x}\phi(x) \| \leq C,
\end{equation}
the associated regularizer satisfies:
\begin{equation}
\|(\nabla^{\ell} R u)(x)\| \leq C' \| \nabla^{\ell} u\|_{\rmL^{1}(\calV^\delta(x))}.
\end{equation}
\end{proposition}
\begin{proof}
By the Deny-Lions lemma there is $C >0$ such that for all $ u \in \rmW^{\ell, 1}(\calV^\delta(x))$.
\begin{equation}
\inf_{f \in \bbP^{\ell -1}} \| u - f \|_{\rmW^{\ell, 1}(\calV^\delta(x))} \cleq   \| \nabla^{\ell} u\|_{\rmL^{1}(\calV^\delta(x))}.
 \end{equation}
The regularizer $R$ preserves the space $\bbP^{\ell -1}$ of polynomials of degree up to $\ell-1$. Combining this with Proposition \ref{prop:nablastab} gives the claimed result.
\end{proof}

The regularizer is adapted to the mesh $\calT_n$ as follows. We choose $\phi_n$ such that for all $x \in T \in \calT_n$:
\begin{align}\label{eq:hypphin}
\phi_n(x) & \ceq h_T,\\
\| \rmD^{1+r} \phi_n(x) \| & \cleq h_T^{-r}, \textrm{ for } 0\leq r \leq \ell.\label{eq:hypphinlast}
\end{align}
We introduce a parameter $\epsilon>0$ and consider the regularizations $R_n^\epsilon = R[\Phi]$ associated with the maps:
\begin{equation}
\Phi_y(x) = x + \epsilon \phi_n(x)y.
\end{equation}
We define:
\begin{equation}
\calV_n^\epsilon(T) = \cup \{ \bbB(x, \epsilon \phi_n(x)) \ : \    x \in T \}.
\end{equation}
We choose $\epsilon$ fixed but small.

\begin{proposition} Fix $\epsilon >0$. For any $T\in \calT_n$ we have an estimate:
\begin{align}
h_T^{1 +(d - \dim T)/q} \| \nabla R_n^\epsilon u\|_{\rmL^q(T)} & \cleq \|u\|_{\rmL^q(\calV_n^\epsilon(T))}, \label{eq:estone}\\
h_T^{(d - \dim T)/q} \| \nabla^\ell R_n^\epsilon u\|_{\rmL^q(T)} & \cleq \|\nabla^\ell u\|_{\rmL^q(\calV_n^\epsilon(T))}. \label{eq:esttwo}
\end{align}
For $T$ of maximal dimension we also have:
\begin{equation}
\| u-  R_n^\epsilon u\|_{\rmL^q(T)}  \cleq h_T^\ell \|\nabla^\ell u\|_{\rmL^q(\calV_n^\epsilon(T))}. \label{eq:estthree}
\end{equation}
\end{proposition}
\begin{proof}
Pick $T \in \calT_n$. Its diameter is $h_T$ and its barycentre $b_T$. Consider the scaling map $\sigma_T: x \to h_T x + b_T$ and let $\hat T$ be the preimage of $T$ by $\sigma_T$, called the reference cell.

Remark that, quite generally, the regularization $R$ transforms as follows under pullback by a diffeomorphism $\sigma$:
\begin{equation}
 R[\Phi] = (\sigma^\star)^{-1} R [\sigma^{-1} \circ \Phi_\bs \circ \sigma ] \sigma^\star,
\end{equation}
where:
\begin{equation}
\sigma^{-1} \circ \Phi_\bs \circ \sigma \ : \ (y,x) \mapsto \sigma^{-1} (\Phi_y ( \sigma(x))).
\end{equation}

In our case the operator $R_n^\epsilon [\sigma_T^{-1} \circ \Phi_\bs \circ \sigma_T ]$ is regularizing on the reference cell $\hat T$ and we denote it by $\hat R_n^\epsilon$. We have:
\begin{equation}
(\sigma_T^{-1} \circ \Phi_y \circ \sigma_T)(x) = x + \epsilon h_T^{-1}\phi_n (h_T x + b_T) y.
\end{equation}
The conditions (\ref{eq:hypphin}), (\ref{eq:hypphinlast}) put us in a position to conclude from Propositions \ref{prop:nablareg} and \ref{prop:nablastab}, that:
\begin{align}
\| \nabla \hat R_n^\epsilon u\|_{\rmL^\infty(\hat T)} & \cleq \|u\|_{\rmL^1(\sigma_T^{-1} \calV_n^{\epsilon} (T))},\\
\| \nabla^\ell \hat R_n^\epsilon u\|_{\rmL^\infty(\hat T)} & \cleq \|\nabla^\ell u\|_{\rmL^1(\sigma_T^{-1}\calV_n^{\epsilon} (T))}.
\end{align}
From this we deduce:
\begin{align}
\| \nabla \hat R_n^\epsilon u\|_{\rmL^q(\hat T)} & \cleq \|u\|_{\rmL^q(\sigma_T^{-1} \calV_n^\epsilon(T))},\\
\| \nabla^\ell \hat R_n^\epsilon u\|_{\rmL^q(\hat T)} & \cleq \|\nabla^\ell u\|_{\rmL^q(\sigma_T^{-1} \calV_n^\epsilon(T))}.
\end{align}
Then the estimates (\ref{eq:estone}) and (\ref{eq:esttwo}) follow from scaling.

Let $T\in \calT_n$ have dimension $d$. From Proposition \ref{prop:regbound} one gets:
\begin{equation}
\| \hat R u \|_{\rmL^q(\hat T)} \cleq \| u\|_{\rmL^{q}(\sigma_T^{-1}\calV_n^{\epsilon} (T)) }.
\end{equation}
Preservation of polynomials and the Deny-Lions lemma then give:
\begin{equation}
\| u - \hat R u \|_{\rmL^q(\hat T)} \cleq \| \nabla^\ell u\|_{\rmL^{q}(\sigma_T^{-1}\calV_n^{\epsilon} (T)) }.
\end{equation}
Scaling then gives (\ref{eq:estthree}).
\end{proof}

For a cell $T \in \calT$ we denote by $\calM_n(T)$ the macroelement surrounding $T$ in $\calT_n$, that is the union of the cells $T' \in \calT_n$ touching $T$:
\begin{equation}
\calM_n(T) = \cup \{ T' \in \calT_n \ : \ T' \cap T \neq \emptyset \}.
\end{equation}

Choose $\epsilon$ so small that for all $n$ and all $T \in \calT_n$,
\begin{equation}
\calV_n^\epsilon(T) \cap S \subseteq \calM_n(T).
\end{equation}

\begin{proposition}
For $u$ defined on $\calV_n^\epsilon(S)$ we have estimates:
\begin{equation}
\| I_n R_n^\epsilon u \|_{\rmL^q(S)} \cleq \|u\|_{\rmL^q(\calV_n^\epsilon(S)}
\end{equation}
and, for $\ell \leq p+1$:
\begin{equation}
\| u - I_n R_n^\epsilon u \|_{\rmL^q(S)} \cleq h_n^\ell \|u\|_{\rmL^q(\calV_n^\epsilon(S)}.
\end{equation}
\end{proposition}

\paragraph{Extension}
We shall define extension operators which extend differential forms on $S$ to some neighborhood of $S$, preserve polynomials up to a certain degree $p$, commute with the exterior derivative and are continuous in $\rmW^{\ell,q}$ norms for $\ell \leq p+1$.

For this purpose we will use maps $\Phi_s$ depending on a parameter $s$, defined outside $S$ with values in $S$ and pull back by these maps. By taking judiciously chosen linear combinations of such pull-backs we meet the requirements of continuity and polynomial preservation.

First we derive some formulas for the derivative of order $\ell$ of the pullback of a differential form. Antisymmetry in the variables will not be important for these considerations, so we consider multilinear rather than differential forms.

Consider then a smooth map $\Phi: \bbR^d \to \bbR^d$ and $u$, a $k$-multilinear form on $\bbR^d$, with $k \geq 1$. We want to give an expression for:
\begin{equation}
(\nabla^\ell \Phi^\star u)[x](\xi_1, \ldots, \xi_{k + \ell}),
\end{equation}
as a linear combination of terms of the form:
\begin{equation}\label{eq:nablaru}
(\nabla^r u)[\Phi(x)]\big(\rmD^{m_1}\Phi(x)(\zeta_1, \ldots, \zeta_{m_1}), \rmD^{m_2}\Phi(x)(\zeta_{m_1 + 1}, \dots, \zeta_{m_1 + m_2}), \dots\big),
\end{equation}
where $m_1 + \ldots + m_{k+r} = k+ \ell$ and:
\begin{equation}
 (\zeta_1, \ldots, \zeta_{k + \ell}) = (\xi_{\sigma(1)}, \ldots, \xi_{\sigma(k + \ell)}),
\end{equation}
for a permutation $\sigma$ of the indexes $\{1, \ldots, k+ \ell \}$.

To handle the combinatorics behind this problem we let $\calD$ be the set of pairs $(\sigma, m)$ where $m= (m_1, \ldots, m_{v(m)})$ is a multi-index of valence $v(m)\geq 1$ and weight $w(m) = m_1 + \ldots + m_{v(p)}$, such that $m_i \geq 1$ for each $1 \leq i \leq v(m)$ ; and $\sigma$ is a permutation of the indexes $(1, \ldots, w(m))$. We define the partial weights:
\begin{equation}
|m|_0= 0 \textrm{ and } |m|_i = m_1 + \ldots + m_i \textrm{ for } 1\leq i \leq v(m).
\end{equation}
We define $(\sigma, m)^\circ[\Phi,x](\xi_1, \ldots, \xi_{w(m)})$ to be the $v(m)$-multivector:
\begin{equation}
\otimes_{i=1}^{v(m)} \rmD^{m_i}\Phi(x)(\zeta_{|m|_{i-1} + 1}, \dots, \zeta_{|m|_i}),
\end{equation}
with:
\begin{equation}
 (\zeta_1, \ldots, \zeta_{w(m)}) = (\xi_{\sigma(1)}, \ldots, \xi_{\sigma(w(m))}).
\end{equation}
Thus:
\begin{equation}\label{eq:sigmacirc}
(\sigma, m)^\circ[\Phi,x] \in \rmL(\otimes^{w(m)} \bbR^d, \otimes^{v(m)}\bbR^d).
\end{equation}
Let $\overline \calD$ be the free group generated by $\calD$. We define some operations in $\overline \calD$, which correspond to differentiating (\ref{eq:sigmacirc}) with respect to $x$.

For a given $(\sigma, m)$ we define, for $j = 1, \ldots, v(m)$, $\partial_j (\sigma,m)$ to consist of the multi-index $(m_1, \ldots, m_j + 1, \ldots, m_{v(m)})$ and the permutation:
\begin{align}
1\leq i \leq w(m) & \mapsto \left\{ \begin{array}{ll}
\sigma (i) &   \textrm{ if }  \sigma(i) \leq |m|_j,\\
\sigma (i) + 1 & \textrm{ if }  \sigma(i) > |m|_j,
\end{array} \right.\\
w(m) + 1 & \mapsto  |m|_j +1.
\end{align}
Defining:
\begin{equation}
\rmD (\sigma,m) = \sum_{j=1}^{v(m)} \partial_j (\sigma, m),
\end{equation}
we have:
\begin{equation}
\rmD ((\sigma,m)^\circ[\Phi,x]) = (\rmD(\sigma,m))^\circ[\Phi,x].
\end{equation}
On the left, $\rmD$ is ordinary differentiation with respect to $x$, and on the right $\rmD$ is an operation in $\overline \calD$.

If $u$ is a $k$-multilinear form and $\zeta = \zeta_1 \otimes \ldots \otimes \zeta_{k} \in \otimes^{k}\bbR^d$ we define the contraction:
\begin{equation}
u : \zeta = u(\zeta_1,\ldots,\zeta_{k}).
\end{equation}
Contraction is bilinear. Given a $k$-multilinear form $u$, we want to differentiate, with respect to $x$,  expressions of the form:
\begin{equation}
(\nabla^r u)[\Phi(x)] : (\sigma, m)^\circ[\Phi,x],
\end{equation}
where $v(m)= k+r$. This corresponds to (\ref{eq:nablaru}).
Since contraction is bilinear and we know how to differentiate the right hand side, it remains to differentiate the left hand side. For this purpose we introduce one more operation in $\overline \calD$. Define $e$ such that $e(\sigma,m)$ consists of the multi-index $(m_1, \ldots, m_{v(m)}, 1)$ and the permutation:
\begin{align}
1\leq i \leq w(m) & \mapsto \sigma(i),\\
w(m)+1 & \mapsto w(m)+1.
\end{align}
Then we have:
\begin{align}\label{eq:nablanabla}
& \nabla \big( (\nabla^r u)[\Phi(x)]  : (\sigma,m)^\circ[\Phi,x]\big)\\
 =&  (\nabla^{r+1} u)[\Phi(x)] : (e(\sigma,m))^\circ[\Phi,x]  +   (\nabla^{r}u)[\Phi(x)] : (\rmD (\sigma,m))^\circ[\Phi,x].
\end{align}

In the free group $\overline \calD$ we now define $\Gamma^\ell_r[k]$ for $\ell \geq 0$ and $0 \leq r \leq \ell$ recursively. We initialize by:
\begin{equation}
\Gamma^0_0[k] = (\id, (1,\ldots,1)),
\end{equation}
with the identity permutation and $k$ terms in the multi-index. Then we define for $\ell \geq 0$:
\begin{align}
\Gamma^{\ell +1}_0[k] & = e \Gamma^{\ell}_0[k],\\
\Gamma^{\ell +1}_{i}[k] & = e \Gamma^{\ell}_i[k] + \rmD \Gamma^{\ell}_{i-1}[k],  \textrm{ for } 1\leq i \leq \ell,\\
\Gamma^{\ell +1}_{\ell +1}[k] & = \rmD \Gamma^{\ell}_{\ell}[k].
\end{align}
One checks that:
\begin{align}
\Gamma^\ell_0[k] & = \Gamma^0_0[k+\ell],\\
\Gamma^\ell_\ell[k] & = \rmD^\ell \Gamma^0_0[k].
\end{align}
\begin{proposition} We have, for a given $k$-form $u$:
\begin{equation}
\nabla^\ell (\Phi^\star u)[x] = \sum_{i= 0}^\ell (\nabla^{\ell-i} u)[\Phi(x)] :\Gamma^\ell_i[k]^\circ[\Phi,x].
\end{equation}
The expression $ \Gamma^\ell_i[k] \in \overline \calD$ is a sum of terms with valence $k+ \ell -i$  and weight $k+ \ell$.
\end{proposition}
\begin{proof}
By induction on $\ell$, using (\ref{eq:nablanabla}).
\end{proof}
For $x$ outside $S$, let $\delta(x)$ denote the distance from $x$ to $S$.

Given a polynomial degree $p$, we will construct an extension operator $E$ as a linear combination of pullbacks:
\begin{align}
E & = \sum_{s \in I} \psi_s \Phi_s^\star, \textrm{ with }
\Phi_s : \calV^\epsilon(S) \setminus S \to S,
\end{align}
subject to the following conditions:
\begin{itemize}
\item The index set $I$ is a finite subset of the interval $[2,3]$ and the coefficients $(\psi_s)_{s\in I}$ are chosen such that for any polynomial $f$ of degree at most $p+ d$:
\begin{equation}
\sum_{s \in I} \psi_s f(s) = f(0).
\end{equation}
\item There is a function $\phi: \calV^\epsilon(S) \setminus S \to \bbR^d$ such that for all $s$ and $x$:
\begin{equation}\label{eq:formofPhi}
\Phi_s (x)  = x + s \phi(x),
\end{equation}
and moreover:
\begin{align}
\|\phi(x)\| & \ceq \delta(x),\\
\| \rmD^{1+ r}\phi(x)\| & \cleq \delta(x)^{-r}, \textrm{ for } 0\leq r \leq p.
\end{align}

\item Finally for $s \in [2, 3]$ we should have:
\begin{align}
\| \rmD \Phi_s(x) ^{-1}\| & \cleq 1,
\end{align}
and the $\Phi_s$ should determine diffeomorphisms from $\calV^\epsilon(S) \setminus S$ for some $\epsilon >0$,  to an interior neighborhood of $\partial S$. 
\end{itemize}

In order to show that the above list of conditions can be met we need some results of a general nature that are given in the Appendix.

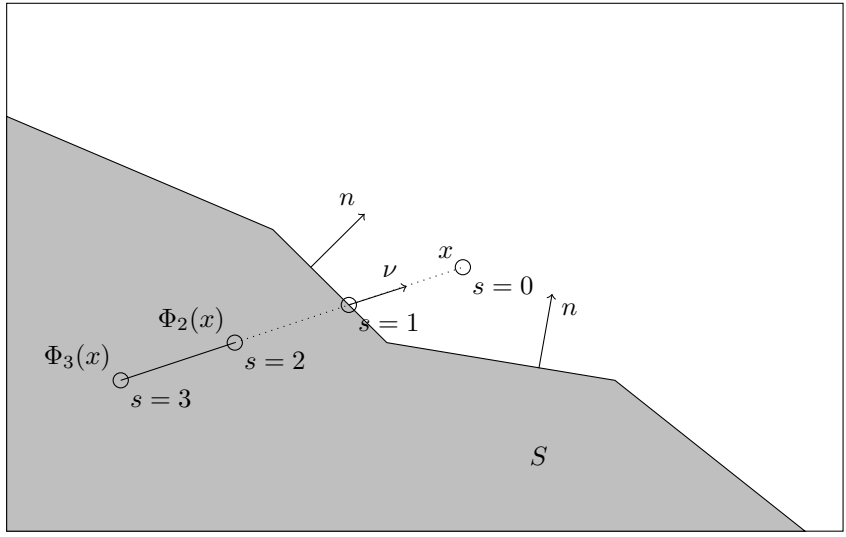
\begin{figure}[htbp]
\begin{center}
\begin{tikzpicture}

\clip[draw] (-1, 0) rectangle (10, 7);

\filldraw[fill= lightgray, draw = black]  (-1, 5.5) -- (2.5, 4) -- (4, 2.5) -- (7, 2) -- (9.5, 0) -- (-1, 0) -- cycle;

\draw[->] (3, 3.5) -- +(45:1cm) node[anchor = south east] {$n$};
\draw[->] (6, 2.16) -- +(80:1cm) node[anchor = north west] {$n$};

\draw[->] (3.5, 3) -- + (18:0.8cm) node[anchor = south east] {$\nu$};

\draw (0.5, 2)   -- (2, 2.5) ;

\draw[dotted] (2, 2.5) -- (3.5, 3)  -- (5, 3.5) ;

\draw (0.5, 2) circle (0.1cm) node[anchor = north west]{$s=3$} node[anchor=south east]{$\Phi_3(x)$};
\draw (2, 2.5) circle (0.1cm) node[anchor = north west]{$s=2$} node[anchor=south east]{$\Phi_2(x)$};
\draw (3.5, 3) circle (0.1cm) node[anchor = north west]{$s=1$};
\draw (5, 3.5) circle (0.1cm) node[anchor = south east]{$x$}  node[anchor = north west] {$s=0$};

\path (6cm,1cm) node{$S$};

\end{tikzpicture}
\end{center}
\caption{Definition of $\Phi_s(x)$, without smoothing. \label{fig:refl}}
\end{figure}

\begin{proposition} The above listed conditions can be met.
\end{proposition}
\begin{proof}
Let $I$ consist of $p+d+1$ points in the interval $[2,3]$. Numbers $\psi_s$ are then determined by solving a linear system with a Vandermonde matrix.

If $\partial S$ had been smooth, the orthogonal projection $\wp$ onto  $\partial S$ would be well defined and smooth on a neighborhood of it. Then we could have taken $\phi(x) = \wp(x) -x$. In the following we modify this construction to allow for Lipshitz boundaries.

Choose a smooth vector field $\nu$ according to Proposition \ref{prop:globalgraph}, pointing outwards on $\partial S$,  so that for some $\epsilon' >0$ the following map is a Lipschitz isomorphism onto its open range:
\begin{equation}
g:\func{\partial S \times ]-\epsilon', \epsilon'[}{\bbR^d}{(z,t)}{ z + t \nu(z)}
\end{equation}
Define $f$ on $\calV^\epsilon(S) \setminus S$ by:
\begin{equation}
f(g(z,t)) = z - g(z,t) = -t \nu(z).
\end{equation}
The problem with $f$, to serve as $\phi$ in (\ref{eq:formofPhi}), is its lack of regularity. 

From Theorem 2 p. 171 in \cite{Ste70} we get a regularized distance function $\tilde \delta$ defined outside $S$, such that:
\begin{equation}
\tilde \delta(x) \ceq \delta (x),
\end{equation}
and for all $r\geq 1$:
\begin{equation}
\| \rmD^r \tilde \delta (x) \| \cleq \delta(x)^{1-r}.
\end{equation}

We regularize $f$ by a variant of (\ref{eq:defreg}, \ref{eq:formphi}). We put:
\begin{equation}
\phi(x) = \int \psi(y) f( x -  \epsilon \tilde \delta(x) y) \rmd y,
\end{equation}
where the parameter $\epsilon$ is chosen to satisfy $\epsilon \tilde \delta(x) \leq 1/2 \delta(x)$, so that $f$ is evaluated far enough from $S$.

For an illustration we refer to Figure \ref{fig:refl}.

The conditions are then met.
\end{proof}

In the following we choose an integer $\ell \leq p$.

\begin{proposition} We have:
\begin{align}
\nabla^\ell E u = &  \nabla^{\ell} u[\Phi_2(x)] \phantom{\sum^\ell_s } \label{eq:lineone}\\
& + \sum_s \psi_s (\nabla^{\ell} u[\Phi_s(x)] - \nabla^\ell u [\Phi_2(x)]): \otimes^{k+\ell}  (\id + s \rmD\phi(x))  \phantom{\sum_s^\ell} \label{eq:linetwo}\\
& + \sum_s \psi_s \int_{0}^{1}\nabla^\ell u [t\Phi_s(x) + (1-t)\Phi_2(x)]:  \phantom{\sum_s^\ell} \label{eq:linethree}\\
& \phantom{+} \qquad \sum_{i=1}^{\ell} \otimes^i(s-2)\phi(x) \otimes \Gamma^\ell_i[k]^\circ[\Phi_s,x] \frac{(1-t)^{i-1}}{(i-1)!}  \rmd t. \label{eq:linefour}
\end{align}
\end{proposition}
\begin{proof}
We have:
\begin{equation}\label{eq:sumnabla}
\nabla^\ell (\Phi_s^\star u)[x] = \sum_{i= 0}^\ell (\nabla^{\ell-i} u)[\Phi_s(x)] :\Gamma^\ell_i[k]^\circ[\Phi_s,x].
\end{equation}
Concerning $\Gamma^\ell_i[k]^\circ[\Phi_s,x]$ remark that:
\begin{align}
\rmD \Phi_s(s) & = \id + s \rmD \phi(x),\\
\rmD^r \Phi_s(s) & =  s \rmD^r \phi(x) \textrm{ for } r \geq 2.
\end{align}
Therefore $\Gamma^\ell_i[k]^\circ[\Phi_s,x]$ is polynomial in $s$ of degree at most  $k + \ell - i $. For $i=0$ the value in $s=0$ is:
\begin{equation}
\Gamma^\ell_0[k]^\circ[\Phi_s,x]|_{s=0} = \otimes^{k+\ell}  (\id + s \rmD\phi(x))|_{s=0} = \otimes^{k+\ell} \id,
\end{equation}
For $i\geq 1$ we have a sum of terms which are products where at least one derivative of order at least $2$ appears, so the value at $s=0$ is $0$.

In (\ref{eq:sumnabla}) the $i=0$ term gives rise to:
\begin{align}
& \sum_s \psi_s \nabla^{\ell} u[\Phi_s(x)]: \Gamma^\ell_0[k]^\circ[\Phi_s,x] =  \nabla^{\ell} u[\Phi_2(x)] +\\
& \qquad \sum_s \psi_s (\nabla^{\ell} u[\Phi_s(x)] - \nabla^\ell u [\Phi_2(x)]): \otimes^{k+\ell}  (\id + s \rmD\phi(x)).
\end{align}
This corresponds to (\ref{eq:lineone}, \ref{eq:linetwo}).

For $i \geq 1 $, Taylor's formula with integral remainder gives:
\begin{align}
& (\nabla^{\ell-i} u)[\Phi_s(x)] = \sum_{j= 0}^{i-1} \nabla^{\ell -i + j}u[\Phi_2(x)]: \otimes^j(s-2)\phi(x) \ \frac{1}{j!} +  \label{eq:sum}\\
& \qquad \int_{0}^{1}\nabla^\ell u [t\Phi_s(x) + (1-t)\Phi_2(x)] : \otimes^i(s-2)\phi(x) \ \frac{  (1-t)^{i-1}}{(i-1)!} \rmd t.
\end{align}
When this expression is  contracted with $\Gamma^\ell_i[k]^\circ[\Phi_s,x]$, the sum (\ref{eq:sum}) consists of polynomials in $s$ with degree $j + k + \ell - i \leq \ell +d $, with value $0$ at $s=0$. Therefore:
\begin{align}
& \sum_s \psi_s (\nabla^{\ell-i} u)[\Phi_s(x)] : \Gamma^\ell_i[k]^\circ[\Phi_s,x]\\
=& \sum_s \psi_s \int_{0}^{1}\nabla^\ell u [t\Phi_s(x) + (1-t)\Phi_2(x)]:\\
& \qquad \otimes^i(s-2)\phi(x) \otimes \Gamma^\ell_i[k]^\circ[\Phi_s,x] \frac{(1-t)^{i-1}}{(i-1)!}  \rmd t.
\end{align}
Summing over $1 \leq i \leq \ell$, this corresponds to (\ref{eq:linethree}, \ref{eq:linefour}).
\end{proof}

From this formula several conclusions can be drawn:
\begin{proposition} \begin{itemize}
\item If $u$ is of class $\calC^\ell(S)$, $Eu$ is of class $\calC^\ell(\calV^\epsilon(S))$.
\item If $u$ is a polynomial of degree at most $p$, $Eu$ also. More precisely if $T$ is a cell touching $\calS$, and $u$ is polynomial on its macro element, of degree at most $p$, then $Eu$ is polynomial of degree at most $p$ on $\calV^\epsilon(T)$.
\item $E$ is bounded  $\rmW^{\ell,q}(S) \to \rmW^{\ell,q}(\calV^\epsilon(S))$. For cells $T\in \calT_n$ touching $\partial S$:
\begin{equation}
\| \nabla^\ell E u \|_{\rmL^{q}(\calV^\epsilon_n(T))} \cleq \| \nabla^\ell u\|_{\rmL^{q}(\calM_n(T))} .
\end{equation} 
\end{itemize}
\end{proposition}
\begin{proof}
For $(\sigma,m) \in \calD$:
\begin{equation}
\| (\sigma, m)^\circ[\Phi_s,x] \, \| \cleq \delta(x)^{v(m) - w(m)},
\end{equation}
from which it follows that:
\begin{equation}
\| \Gamma^\ell_i[k]^\circ[\Phi_s,x]\, \|  \cleq \delta(x)^{-i}.
\end{equation}
This gives:
\begin{equation}
\| \otimes^i\phi(x) \otimes \Gamma^\ell_i[k]^\circ[\Phi_s,x]\, \| \cleq 1. \label{eq:gammabound}
\end{equation}
Concerning continuity properties, we check it on $\partial S$. Suppose $x_0 \in \partial S$ and $x \to x_0$. Then we have:
\begin{equation}
\| \nabla^{\ell} u[\Phi_s(x)] - \nabla^\ell u [\Phi_2(x)] \, \| \to 0,
\end{equation}
so that (\ref{eq:linetwo}) converges to $0$. If the integrals (\ref{eq:linethree}) had been evaluated at $\nabla^\ell u[\Phi_2(x)]$ their sum over $s\in I $ would be $0$. A similar argument to the above, combined with (\ref{eq:gammabound}), then shows that (\ref{eq:linethree}, \ref{eq:linefour}) also converges to $0$. We are left with the term on (\ref{eq:lineone}) which converges to $\nabla^\ell u (x_0)$.

Boundedness properties of $E$ follow from (\ref{eq:gammabound}) and the assumption that the $\Phi_s$ determine diffeomorphisms with uniformly bounded Jacobian determinants.
\end{proof}

\paragraph{Quasi-interpolator}
Putting together the pieces we get:
\begin{theorem} For any $\epsilon > 0$ the operators $Q^\epsilon_n = I_nR_n^\epsilon E$ satisfy local estimates, for $T \in \calT^d_n$:
\begin{align}
\|Q^\epsilon_n u \|_{\rmL^q(T)} & \cleq \| u \|_{\rmL^q(\calM_n(T))},\\
\| u - Q^\epsilon_n u \|_{\rmL^q(T)} & \cleq h_T^\ell \| \nabla^\ell u \|_{\rmL^q(\calM_n(T))},
\end{align}
as the corresponding global ones:
\begin{align}
\|Q^\epsilon_n u \|_{\rmL^q(S)} & \cleq \| u \|_{\rmL^q(S)},\\
\| u - Q^\epsilon_n u \|_{\rmL^q(S)} &\cleq h^\ell \| \nabla^\ell u \|_{\rmL^q(S)}.
\end{align}
Moreover for any $\epsilon'$, choosing $\epsilon$ small enough will yield, for $u \in X^k_n$:
\begin{align}
\| u - Q^\epsilon_n u \|_{\rmL^q(T)} &\leq \epsilon'  \| u \|_{\rmL^q(\calM_n(T))},\\
\| u - Q^\epsilon_n u \|_{\rmL^q(S)} &\cleq \epsilon'  \| u \|_{\rmL^q(S)}.
\end{align}
Finally $Q^\epsilon_n$ commutes with the exterior derivative (when it is in $\rmL^q(S)$). 
\end{theorem}

When $\epsilon$ is chosen so small that $\| (\id - Q^\epsilon_n)|_{X^k_n} \|_{\rmL^q(S) \to \rmL^q(S)} \leq 1/2$,  $Q^\epsilon_n|_{X^k_n}$ is invertible with norm less than $2$. We define operators:
\begin{equation}
P_n = (Q^\epsilon_n|_{X^k_n})^{-1} Q^\epsilon_n.
\end{equation}

\begin{proposition} The operators $P_n$ are $\rmL^q(S)$ stable projections onto $X^k_n$ which commute with the exterior derivative.
\end{proposition}

The case $q = 2$ leads to eigenvalue convergence for the operator $\rmd^\star \rmd$ discretized by the Galerkin method on $X^k_n$ and therefore for the Hodge-Laplacian in mixed form. For a discussion of eigenvalue convergence we refer to \cite{Bof10}, \cite{ArnFalWin10} and \cite{ChrWin13IMA}.


\section{Sobolev injection and translation estimate \label{sec:sob}}
In this section we prove a Sobolev injection theorem generalizing the one we introduced in \cite{ChrSch11}. The technique of proof is slightly different, and we generalize to differential forms in all dimensions. We also prove a translation estimate of the type introduced in \cite{KarKar11}, Theorem A.1. Compared with that paper we get an optimal bound and a generalization to all known mixed finite elements in the $h$-version.

Given a cellular complex $\calT$ we define a broken $\rmH^1$ seminorm as follows:
\begin{equation}
\llp u \rrp^2 = \sum_{T \in \calT^d} \| \nabla u \|^2_{\rmL^2(T)} + \sum_{T \in \calT^{d-1}} h_T^{-1} \| [u]_T \|^2_{\rmL^2(T)} .
\end{equation}
Given a simplex $T \in \calT^{d-1}$, $[u]_T$ denotes the jump of $u$ across $T$. The scaling factor in front of the jump terms is chosen so that the two terms that are summed scale in the same way, see Lemma \ref{lem:scaling}.

We suppose that we have a sequence of cellular complexes $(\calT_n)$ and that each $\calT_n$ is equipped with a compatible FE system $A[n]$. We suppose that they are of the type discussed in the previous section so that those constructions apply. We define $X^k_n= A[n]^k(\calT_n)$. The above defined broken seminorm will be relative to one of the $\calT_n$ and it will be clear from the context which one, so that we omit the precision from the notation.

In the $\rmL^2(S)$ case define:
\begin{align}
X^k & = \{ u \in \rmL^2(S) \ : \ \rmd u \in \rmL^2(S) \},\\
W^k & = \{ u \in X^k \ : \ \rmd u = 0 \},\\
V^k & = \{ u \in X^k \ : \ \forall w \in W^k \quad \ts \int u \cdot w = 0 \}.
\end{align}

Define:
\begin{align}
W^k_n & = \{ u \in X^k_n \ : \ \rmd u = 0 \},\\
V^k_n & = \{ u \in X^k_n \ : \ \forall w \in W^k_n \quad \ts \int u \cdot w = 0 \}.
\end{align}

Denote by $\hodge$ the $\rmL^2$ orthogonal projection onto $\overline V^k$, the completion of $V^k$ in $\rmL^2(S)$. This operator realizes a Hodge decomposition of $u$ in the form $u = (u - \hodge u) + \hodge u$. The following well-known trick is also useful in the proof of eigenvalue convergence.
\begin{proposition} We have, for $u \in V^k_n$:
\begin{equation}
\| u -\hodge u \|_{\rmL^2(S)} \leq \| \hodge u - P_n \hodge u\|_{\rmL^2(S)}.
\end{equation}
\end{proposition}
\begin{proof}
We have:
\begin{equation}
u - P_n \hodge u = P_n(u -\hodge u) \in W^k_n.
\end{equation}
Now write:
\begin{equation}
\hodge u - P_n \hodge u = (u - P_n \hodge u) - (u - \hodge u),
\end{equation}
and remark that the two terms on the left hand side are orthogonal.
\end{proof}

\begin{proposition} For $u \in X^k_n$:
\begin{align}
\llp u \rrp & \ceq \| \nabla R^\epsilon_n E u \|_{\rmL^2(S)},\label{eq:llpequiv}\\
\| u - R^\epsilon_n E u \|_{\rmL^2(S)} & \cleq h \llp u \rrp.
\end{align}
\end{proposition}
\begin{proof}
The parameter $\epsilon$ is chosen so small that for all $u \in X^k_n$:
\begin{equation}
\llp u - R^\epsilon_n E u \rrp \leq 1/2 \llp u \rrp.
\end{equation}
Then we have:
\begin{align}
\llp u \rrp & \leq 2 \llp R^\epsilon_n E u \rrp,\\
\llp R^\epsilon_n E u \rrp & \leq 3/2 \llp u \rrp.
\end{align}
Since $R^\epsilon_n E u$ is smooth we have:
\begin{equation}
\llp R^\epsilon_n E u \rrp = \| \nabla R^\epsilon_n E u \|_{\rmL^2(S)}.
\end{equation}
This gives (\ref{eq:llpequiv}).

The other estimate is proved locally by scaling from a reference macroelement.
\end{proof}

\begin{proposition} For $u \in \rmH^1(S)$:
\begin{equation}
\llp P_n u \rrp \cleq \| \nabla u \|_{\rmL^2(S)}.
\end{equation}
\end{proposition}
\begin{proof}
By restriction to reference simplexes and scaling:
\begin{equation}
\llp Q_n u \rrp \cleq \| \nabla u \|_{\rmL^2(S)}.
\end{equation}
Choosing $\epsilon$ small enough we also have for $u \in X^k_n$:
\begin{equation}
\llp u - Q_n u \rrp \leq 1/2 \llp u \rrp.
\end{equation}
Combining the two we get the proposition.
\end{proof}

\begin{proposition} Suppose $S$ is convex and that the meshes are quasiuniform. Then for all $u \in V^k_n$:
\begin{equation}\label{eq:onenbd}
\llp u \rrp \cleq \| \rmd u \|_{\rmL^2(S)}. 
\end{equation}
\end{proposition}
\begin{proof}
We have, for $u \in V^k_n$:
\begin{align}
\llp u \rrp & = \llp P_n u \rrp \\
& \leq \llp P_n (u - \hodge u)  \rrp +  \llp P_n \hodge u \rrp, \\
& \leq h^{-1} \| P_n (u - \hodge u) \|_{\rmL^2} + \| \nabla \hodge u \|_{\rmL^2},\\
& \cleq h^{-1} \| u - \hodge u\|_{\rmL^2} + \| \nabla \hodge u \|_{\rmL^2},\\
& \cleq \| \nabla \hodge u \|_{\rmL^2}.
\end{align}
From this (\ref{eq:onenbd}) follows.
\end{proof}

\begin{proposition} Let $q$ be the relevant Sobolev exponent, so that $\rmH^1(S) \to \rmL^q(S)$. For all $u \in X^k_n$:
\begin{equation}
\| u \|_{\rmL^q(S)} \cleq \llp u \rrp + \| u \|_{\rmL^2(S)}.
\end{equation}
\end{proposition}
\begin{proof} We have:
\begin{align}
\| u \|_{\rmL^q(S)} & \cleq \| R^\epsilon_n E u \|_{\rmL^q(S)},\\
&\cleq  \| \nabla R^\epsilon_n E u\|_{\rmL^2(S)} + \| R^\epsilon_n E u\|_{\rmL^2(S)},\\
&\cleq \llp u \rrp + \| u \|_{\rmL^2(S)},
\end{align}
as claimed.
\end{proof}

Denote by $\tau_y$ the translation by the vector $y$, so that when $u$ is defined in $x-y$ we have:
\begin{equation}
(\tau_y u )(x) = u(x -y).
\end{equation}

\begin{proposition} For all $u \in X^k_n$:
\begin{equation}
\| u - \tau_y E u \|_{\rmL^2(S)} \cleq (|y| + h^{1/2} |y|^{1/2}) \llp u \rrp.
\end{equation}
\end{proposition}
\begin{proof}
On a reference simplex $\hat T$ we can write for $u \in X^k_n$ pulled back:
\begin{equation}
\| \hat u - \tau_{\hat y} \hat u \|^2_{\rmL^2(\hat T)} \cleq |\hat y|^2 \sum_{T' \in \calM_n(\hat T)^d} \| \nabla \hat u \|^2_{\rmL^2(T')} + |\hat y| \sum_{T' \in \calM_n(\hat T)^{d-1}} \| [\hat u]_{T'}  \|^2_{\rmL^2(T')}.
\end{equation}
Scaling back to $T$ of size $h$ we get, with $y = h \hat y$:
\begin{equation}
\| u - \tau_{y} u \|^2_{\rmL^2(T)} \cleq |y|^2 \sum_{T' \in \calM_n(T)^d} \| \nabla  u \|^2_{\rmL^2(T')} + |y| \sum_{T' \in \calM_n(T)^{d-1}} \| [ u]_{T'}  \|^2_{\rmL^2(T')},
\end{equation}
so that:
\begin{equation}
\| u - \tau_{y} u \|^2_{\rmL^2(T)} \cleq (|y|^2 + h |y|) \llp u  \rrp^2_{\calM_n(T)}.
\end{equation}
This estimate comes with a restriction of the type $|y| \leq  h / C$ ensuring that one does not translate $T$ out of its associated macro element. For $|y| \geq h/C$ we can write:
\begin{equation}
u - \tau_y E u = (u - R^\epsilon_n E u) + (R^\epsilon_nEu - \tau_y ER^\epsilon_nEu ) +  \tau_yE(R^\epsilon_n Eu -u).
\end{equation}
From this we deduce:
\begin{align}
\| u - \tau_y E u \|_{\rmL^2} & \cleq  \|u - R^\epsilon_n E u\|_{\rmL^2}  + \| R^\epsilon_nEu - \tau_y ER^\epsilon_nEu \|_{\rmL^2},\\
&\cleq  h \llp u \rrp + |y| \| \nabla R^\epsilon_nE u \|_{\rmL^2},\\
&\cleq |y| \, \llp u \rrp.
\end{align}
This concludes the proof.
\end{proof}

\section*{Appendix}

\begin{lemma}\label{lem:lipinv}
In some Banach space $\bbE$, let $U$ be an open set and $f:U \to \bbE$ be a short map, i.e. for some $\delta <1$:
\begin{equation}\label{eq:short}
\|f(x) - f(y)\| \leq \delta \|x-y\|.
\end{equation}
Then the map $g:U \to \bbE$, $ x \mapsto x + f(x)$ has an open range $V$ and determines a Lipschitz bijection $U \to V$, with a Lipschitz inverse.
\end{lemma}
\begin{proof}
Suppose that $g(x_0) = y_0$. For $\|y - y_0\| \leq \epsilon$ find the solution $x$ of $g(x) = y$ as a fixed point of the map $z \mapsto y - f(z)$. More precisely construct a sequence starting at $x_0$ and defined by $x_{n+1} = y - f(x_n)$. Then, as long as it is defined (in $U$) we have: 
\begin{align}
\|x_{n+1} - x_n\| &\leq \delta \| x_n - x_{n-1} \|,\\
&\leq \delta^{n-1} \epsilon.
\end{align}
If $\epsilon$ is chosen so small that the closed ball with centre $x_0$ and radius $(1-\delta)^{-1} \epsilon$ is included in $U$, the sequence is defined for all $n$, and converges to a limit $x\in U$ solving $g(x) = y$.

It follows that $g$ is an open mapping. In particular the range $V$ is open. Moreover:
\begin{equation}
\|g(x) - g(y) \| \geq (1- \delta) \|x-y\|.
\end{equation}
This gives injectivity, so that $g: U \to V$ is bijective. The inverse is Lipschitz with constant not worse than $(1-\delta)^{-1}$.
\end{proof}

\begin{lemma}\label{lem:directionchange}
In some Euclidean space $\bbE$, let $S$ be a bounded domain, whose boundary $\partial S$ is locally the graph of a Lipschitz function. Let $n$ be the outward pointing normal on $\partial S$. Suppose $m$ is a unit vector, that $x_0 \in \partial S$ and that for $x$ in a  neighborhood of $x_0$ in $\partial S$ we have $n(x) \cdot m \geq \epsilon$, for some $\epsilon >0$. Then there is a neighborhood of $x_0$ in $\partial S$ which is a Lipschitz graph above the plane orthogonal to $m$.
\end{lemma}

\begin{proof}
We know that for a certain outward pointing unit vector $m_0$, a neighborhood $\calU_0$ of $x_0$ is a Lipschitz graph above an open ball $B_0$ in $m_0^\perp$. Choose $\theta \in [0, \pi/2[$ such that, for $x \in \calU_0$, $n(x) \cdot m_0 \geq \cos \theta$ and moreover $\epsilon \geq \cos \theta$. Let $f : B_0 \to \bbR$ be the function such that $\calU_1$ is the range of:
\begin{equation}
 \func{B_0}{\bbE}{y}{y + f(y)m_0}
\end{equation}
Since, for $y \in B_0$:
\begin{equation}
n(y + f(y)m_0) = (m_0 - \grad f(y))/(1 + |\grad f(y)|^2)^{1/2},
\end{equation}
we get $|\grad f(x)| \leq \tan \theta$.

Choose $y,y' \in B_0$ and put $x= y + f(y)m_0$ and $x'= y' + f(y')m_0$. For $s,s' \in \bbR$ we have:
\begin{align}
&|(x + sm_0) - (x' + s' m_0)|^2\\
 =& |y-y'|^2 + (f(y)-f(y'))^2 + 2(f(y)- f(y'))(s-s') + (s-s')^2.
\end{align}
Then remark that for $M >0$:
\begin{align}
 2(f(y)- f(y'))(s-s') \leq& (1 + M^{-2}) (f(y)-f(y'))^2 + (1 + M^{-2})^{-1}(s-s')^2\\
\leq& (f(y)-f(y'))^2 + M^{-2} \tan^2 \theta |y-y'|^2\\
&+ (1 + M^{-2})^{-1}(s-s')^2. 
\end{align}
In particular with $M = \tan \theta$ we get:
\begin{align}
|(x + sm_0) - (x' + s' m_0)|^2 \geq  (1 - (1 + \tan^{-2} \theta)^{-1}) (s - s')^2,
\end{align}
which simplifies to:
\begin{align}
|(x + sm_0) - (x' + s' m_0)| \geq  \cos \theta |s - s'|,
\end{align}
Define the function:
\begin{equation}\label{eq:gzero}
g_0: \func{\calU_0 \times \bbR}{\bbE}{x}{y + sm_0}.
\end{equation}
Its range is $B_0 + \bbR m_0$, and it is bijective onto it.

Consider now a unit vector $m_1$ such that $n(x) \cdot m_1 \geq \cos \theta$ for $x \in \calU_0$. Define $g_1$ as in (\ref{eq:gzero}), replacing $m_0$ by $m_1$. We shall show that $g_1$ is open and bi-Lipschitz, when $|m_1 - m_0| < \cos \theta$. Define $f_1$ on  $B_0 + \bbR m_0$, by  $f_1(x) = g_1 \circ g_0^{-1} - x$. We have, with the preceding notations:
\begin{equation}
f_1(x + sm_0) - f_1(x' + s'm_0)  = (s-s')(m_1 - m_0),
\end{equation}
hence:
\begin{equation}
|f_1(x + sm_0) - f_1(x' + s'm_0)| \leq |m_1 - m_0| / \cos \theta |(x + sm_0) - (x' + s' m_0)|.
\end{equation}
One then applies Lemma \ref{lem:lipinv} and deduces that $g_1$ is open and bi-Lipschitz. It follows that there is a ball $B_1$ in $m_1^\perp$ above which a neighborhood $\calU_1 \subseteq \calU_0$ of $x_0$ is a graph. We may repeat the above considerations to construct a sequence of such vectors $m_1, m_2, \ldots, m_k$ reaching $m$ in a finite number of steps.
\end{proof}

\begin{proposition}\label{prop:globalgraph}
In some Euclidean space $\bbE$, let $S$ be a bounded domain, whose boundary $\partial S$ is locally the graph of a Lipschitz function. Then there exists a smooth vector field $\nu$ on $\bbE$, of unit length and outward pointing on $\partial S$, such that for some $\epsilon >0$ the map:
\begin{equation}\label{eq:magicgraph}
\func{\partial S \times ]-\epsilon, \epsilon[}{\bbE}{(x,s)}{ x + s \nu(x)}
\end{equation}
has open range and determines a Lipschitz isomorphism onto it.
\end{proposition}
\begin{proof}
Cover $\partial S$ with a finite number of orthogonal cylinders $C_i$ directed along a unit vector $n_i$ pointing out of $S$, and with a base $U_i$, such that above $U_i$, $\partial S$ is the graph of a Lipschitz function. Let $n$ denote the outward pointing normal on $\partial S$. For some $\theta \in [0, \pi/2[$ we have for all $i$ and all $x \in \partial S \cap C_i$: 
\begin{equation}\label{eq:transverse}
n_i \cdot n(x) \geq \cos (\theta).
\end{equation}
Choose smooth functions $\alpha_i$ on $\bbE$ whose restriction to $\partial S$ form a partition of unity. Define:
\begin{equation}
\tilde \nu(x) = \sum_i \alpha_i(x) n_i,
\end{equation}
and normalize by putting:
\begin{equation}
\nu(x) = \tilde \nu(x) / |\nu (x)|.
\end{equation}
We have $\nu(x) \cdot n(x) \geq \cos( \theta)$, for all $x \in \partial S$.

Denote by $g$ the function:
\begin{equation}
g: \func{\partial S\times \bbR }{\bbE}{(x,s)}{ x + s\nu(x) }
\end{equation}

Pick $x_0 \in \partial S$ and put $m_0=\nu(x_0)$. By Lemma \ref{lem:directionchange}, $\partial S$ is locally a Lipschitz graph above the plane orthogonal to $m_0$. Let $\calU(x_0)$ be the corresponding neighborhood of $x_0$ in $\partial S$.
Denote by $g_0$ the function:
\begin{equation}
g_0: \func{\partial S \cap \calU(x_0) \times \bbR }{\bbE}{(x,s)}{ x + sm_0 }
\end{equation}
It is a Lipschitz isomorphism onto its range, which is open in $\bbE$.

Define $f$ by $f(x)= g \circ g_0^{-1}(x) -x$. We have:
\begin{equation}
f(x+ sm_0) = s(\nu(x) - m_0),
\end{equation}
so that:
\begin{align}
f(x+ sm_0) - f(y + tm_0) & = (s-t)(\nu(x) - m_0) + t(\nu(x) - \nu(y)).
\end{align}
It follows that, for a small enough $\epsilon$ and possibly reducing $\calU(x_0)$, $f$ is short on $g_0(\calU(x_0) \times ]-\epsilon, \epsilon[)$. By Lemma \ref{lem:lipinv} it follows that $g \circ  g_0^{-1}$ restricted to  $g_0(\calU(x_0) \times ]-\epsilon, \epsilon[)$, for some small enough $\epsilon$, called $\epsilon(x_0)$, has open range and determines a Lipschitz isomorphism onto it. 

Hence $g$ restricted to $\calU(x_0) \times ]-\epsilon(x_0), \epsilon(x_0)[$ has open range and determines a Lipschitz isomorphism onto it. In particular $g$ is an open mapping.

The open subsets $\calU(x_0)$ associated with each $x_0 \in \partial S$ cover $\partial S$. Choose a finite subset $\calF$ of $\partial S$  such that the sets $\calU(x)$ for $x \in \calF$ cover $\partial S$. Pick $\mu>0$ such that if $|x-x'| \leq \mu$ they belong to a common $\calU(x)$ for $x \in \calF$. In (\ref{eq:magicgraph}), choose $\epsilon$ smaller than each $\epsilon(x)$ for $x \in \calF$, and also smaller than $\mu/3$. Pick now two points $x,x' \in \partial S$, and $s,s'$ in  $]-\epsilon, \epsilon[$. If $|x-x'| \leq \mu$ they belong to a common $\calU(x)$, $x \in \calF$. If on the other hand $|x-x'| \geq \mu$ we have:
\begin{align}
|g(x,s) - g(x',s')| &\geq |x-x'| - |s| - |s'|,\\
&\geq \mu/3.
\end{align}
Based on these two cases we may conclude that for some global $m$:
\begin{equation}
|g(x,s) - g(x',s')| \geq m (|x-x'|^2 + |s-s'|^2).
\end{equation}
The lemma follows.
\end{proof}

\section{Acknowledgments} On finite elements, I have benefited from the insights of in particular Annalisa Buffa, Jean-Claude N\'ed\'elec and Ragnar Winther. On algebraic topology I am grateful for the help of John Rognes and Bj\o rn Jahren, in particular with Propositions \ref{prop:gluecoh} and \ref{prop:tensorext}. I am also grateful to Martin Costabel for a helpful introduction to Stein's construction of universal extension operators. Example \ref{ex:euler} is from a discussion with Thomas Dubos (shallow water on the sphere) during spring 2009. Trygve Karper's remarks on upwinding were also helpful.

This work, conducted as part of the award ``Numerical analysis and simulations of geometric wave equations'' made under the European Heads of Research Councils and European Science Foundation EURYI (European Young Investigator) Awards scheme, was supported by funds from the Participating Organizations of EURYI and the EC Sixth Framework Program.

\bibliography{../../../Bibliography/alexandria,../../../Bibliography/mybibliography}{}
\bibliographystyle{plain}

\end{document}